\documentclass{svjour3}

\usepackage{cite}
\usepackage{url}
\usepackage{amsmath}
\usepackage{amsfonts}
\usepackage{amssymb}
\usepackage{setspace}
\usepackage{float}
\usepackage{graphicx}
\usepackage[dvipsnames]{xcolor}

\usepackage{caption}

\usepackage{algorithm}
\usepackage{algpseudocode}

\newcommand{\post}{\mathbb{P}_{\textit{post}}}
\newcommand{\like}{\mathbb{P}_{\textit{like}}}
\newcommand{\prior}{\mathbb{P}_{\textit{prior}}}

\newcommand{\q}{\textbf{q}}
\newcommand{\pars}{p,z_{0},L}

\newcommand{\x}{\textbf{x}}
\newcommand{\w}{\textbf{w}}
\newcommand{\vv}{\textbf{v}}
\newcommand{\argmax}{\mathop{\mathrm{argmax}}}

\newcommand{\dv}{\nabla\cdot}

\usepackage[colorinlistoftodos,prependcaption,textsize=tiny]{todonotes}

\newcommand{\hll}[1]{{\color{black} #1}}

\newcommand{\mytheta}{\boldsymbol{\theta}}
\newcommand{\pkgname}[1]{{\tt #1}}
\newcommand{\hl}[1]{{#1}}

\usepackage[top=2.5cm, bottom=2.5cm, right=3cm, left=2.2cm]{geometry}

\DeclareGraphicsExtensions{.pdf,.jpg,.png}

\begin{document}

\journalname{Pure and Applied Geophysics}
\date{\today}
\title{Simultaneous model calibration and source inversion in
  atmospheric dispersion models} 

\author{Juan G. Garc\'{i}a \and\ Bamdad Hosseini \and\ John M. Stockie}
\institute{J.G. Garc\'{i}a 
\at 3DM Devices Inc. 26019 31B, Aldergrove, BC V4W 2Z6, Canada\\
\email{juan\_g@3dm.com}
  \and\ J.M. Stockie
  \at Department of Mathematics, Simon Fraser University, 8888
  University Drive, Burnaby, BC V5A 1S6, Canada\\
  \email{jstockie@sfu.ca}
  \and
  B. Hosseini
  \at Computing and Mathematical Sciences, California Institute of
  Technology, 1200 E. California Blvd., Pasadena, CA 91125, USA\\
  \email{bamdadh@caltech.edu}
}

\maketitle

\begin{abstract}
  We present a cost-effective method for model calibration and solution
  of source inversion problems in atmospheric dispersion modelling.  We
  use Gaussian process emulations of atmospheric dispersion models
  within a Bayesian framework for solution of inverse problems. The
  model and source parameters are treated as unknowns and we obtain
  point estimates and approximation of uncertainties for sources while
  simultaneously calibrating the forward model.  The method is validated
  in the context of an industrial case study involving emissions from a
  smelting operation for which cumulative monthly measurements of zinc
  particulate depositions are available.
  \keywords{%
    Atmospheric dispersion \and
    Particle deposition \and
    Inverse source identification \and
    Bayesian estimation \and 
    Uncertainty quantification}
  \subclass{%
    62F15 \and %
    65M08 \and
    65M32 \and %
    86A10}
\end{abstract}

\section{Introduction}

\hl{ This article presents a methodology for simultaneously estimating
  pollutant emission rates and calibrating atmospheric dispersion models
  using far-field measurements of particulate deposition.  We use a
  Gaussian process (GP) emulator to efficiently approximate a relatively
  expensive partial differential equation (PDE) model for atmospheric
  dispersion which is then coupled with a Bayesian framework for source
  inversion and uncertainty quantification.

  Source identification (or inversion) problems are prevalent in various
  areas of environmental monitoring \cite{kim2007advanced} and forensics
  \cite{mudge2008environmental}. These problems are often difficult to
  solve due to the scarcity of data and the complexity of physical
  processes that obscure the true nature of sources. Traditionally,
  source identification problems are solved by a combination of
  empirical techniques such as chemical analysis
  \cite{morrison2000application}, remote sensing methods like thermal
  imaging \cite{matson1981identification}, and statistical tools such as
  principle component analysis \cite{mudge2007multivariate}.  Owing to 
  the prevalence of measurement methods and chemical analysis
  in the forensics literature \cite{colbeck2009, hesterharrison2008,
    spikmans2019forensics}, the application of mathematical models and
  numerical simulations is relatively under-developed. This is in part due
  to complexity of certain mathematical models, uncertainty in values of
  tuning parameters, and general difficulty in assessing confidence in
  simulations.  And even when the errors and uncertainties in contaminant 
  measurements may be well-known or controlled, the effect of those
  errors when they are used as inputs for a source inversion model is
  much less well understood.  It is therefore crucial to quantify
  the effect of measurement errors on emissions estimates and obtain
  uncertainty bounds in order to make reliable and defensible claims
  in forensics applications.

  In the context of atmospheric pollutant dispersion, numerous
  mathematical models are available that provide approximations to} the
actual physical processes driving dispersion of particulates and other
substances within the atmosphere \cite{chatwin1982use,
  lewellen1989meteorological}. A common challenge in dispersion
modelling is to maintain accuracy since many models suffer from sizable
errors and large uncertainties in parameters.  In some cases, these
inaccuracies can be addressed by increasing the model complexity to
account for physical processes that have been neglected, in the hope
that the resulting more complex model yields correspondingly higher
accuracy and more realistic simulations.  However, in many practical
applications, such extensions lead to no significant improvement in the
solution because the model errors are far outweighed by the uncertainty
in values of the model parameters.  This difficulty is further
exacerbated by the fact that many models contain a large number of
unknown model parameters.  Some of these parameters have a clear
physical basis, such as diffusion coefficients, particulate density and
settling velocities, whereas others are mathematical (or fitting)
parameters such as the Monin-Obukhov length, atmospheric stability
classes, or terrain roughness length.

Direct measurements of many physical and mathematical parameters are
often problematic or even infeasible; for example, eddy diffusion
coefficients are typically difficult to estimate in practice
\cite{seinfeld1998atmospheric}. One is therefore often forced to deal
with significant uncertainties in the true value of model parameters,
even when complex models are employed that under ideal circumstances
could potentially simulate results very accurately.  Because large
errors in model parameters can severely contaminate simulations, it is
therefore essential when simulating atmospheric dispersion scenarios to
simultaneously understand and control the parameter uncertainty.  Sources
of uncertainty are diverse and can be classified roughly into one of
three types \cite{rao2005uncertainty}:
\begin{itemize}
\item Data uncertainty: when the empirical measurements or model
  parameter values are inaccurate; 
\item Model uncertainty: where the model does not capture all physical
  processes of interest; or 
\item Stochastic uncertainty: relating to the inherent unpredictability
  or randomness in physical processes, such as atmospheric turbulence.
\end{itemize}
In this work we focus primarily on dealing with the first case of {data
  uncertainty} appearing in the context of atmospheric source inversion
problems.

Despite the difficulty of obtaining accurate simulations, there is a
growing need to develop more accurate and computationally efficient
models for practical atmospheric dispersion applications
\cite{leelHossy2014dispersion}.  The necessity for computational
efficiency comes from the fact that models for atmospheric transport are
usually expressed as systems of partial differential equations that are
computationally expensive to solve at the typical resolutions required,
which is a particular limitation for real-time applications.  Efficiency
is of even greater importance in the context of source inversion
problems, since a single solution of the inverse problem may require a
large number of evaluations of the forward dispersion model.  One of the
goals of this work is therefore to demonstrate how to use
computationally demanding forward models to solve atmospheric source
inversion problems in an efficient manner.

\subsection{Atmospheric dispersion modelling}

We model the dispersion of pollutants in the atmosphere using the
advection-diffusion PDE
\hl{
\begin{gather}
  \label{eqnModelIntro}
  \frac{\partial C(\x,t)}{\partial t} + 
  \dv(\bar{\vv}C(\x,t)-\textbf{D}\nabla C(\x,t)) = f(\x,t), 
\end{gather}
} where $C(\x,t)$ denotes the concentration of pollutant at location
$\x=(x,y,z)\in\mathbb{R}^{3}$, $t$ is time, $\bar{\vv}$ is the wind
velocity field, $\textbf{D}$ is the diffusivity tensor, and $f(\x,t)$ is
a pollutant source term. \hl{Throughout this paper we always consider
  \eqref{eqnModelIntro} with zero initial condition $C(\x, 0) = 0$, so
  that no pre-existing transient concentrations exist. We make this
  assumption for reasons of simplicity, and because in practice it is
  often very difficult to measure these pre-existing concentrations.} We
focus on the case of a finite collection of $n$ point sources each
having an emission rate that is constant in time, which allows $f$ to be
written as
\begin{gather}
  f(\x,t) = \sum_{j=1}^{n} q_{j} \delta(\x-\x_{j}),
\end{gather}
where $\q = (q_1, q_2, \ldots, q_n)^T \in \mathbb{R}^n$ represent
emission rates, \hl{$\{\x_j\}_{j=1}^n$ are the locations of the sources},
and $\delta(\x-\x_{j})$ are Dirac delta functions representing point sources at the
locations $\x_j$.  The wind velocity field
$\bar{\vv}$ and diffusivity tensor $\textbf{D}$ are often difficult to
measure and so are replaced using simpler mathematical approximations
that depend upon empirical parameters describing the variation of
$\bar{\vv}$ and $\textbf{D}$ with time and space. Common examples of
such parameters include roughness length, Monin-Obukhov length and
Pasquill stability class \cite{seinfeld1998atmospheric}.  We combine all
empirical parameters together into a vector $\mytheta = (\theta_k) \in
\mathbb{R}^m$, and then denote the linear operator
$\dv(\vv+\textbf{D}\nabla)$ on the left hand side of
$\eqref{eqnModelIntro}$ by $\textbf{L}(\mytheta)$ to make explicit the
dependence of the operator on the parameter vector $\mytheta$, noting
that the dependence of $\textbf{L}$ on $\mytheta$ can still be
nonlinear.  Then we rewrite \eqref{eqnModelIntro} as
\begin{gather}
  \label{eqnParametric}
  \frac{\partial C(\x,t)}{\partial t} + \textbf{L}(\mytheta)C(\x,t) =
  \sum_{j=1}^{n} q_{j} \delta(\x-\x_{j}). 
\end{gather}
The task of finding the concentration $C(\x,t)$ given $\q$ and
$\mytheta$, supplemented by suitable boundary and initial conditions,
is referred to as the {\it forward problem} in the context of
atmospheric pollutant transport.

\subsection{Atmospheric source inversion}

The goal of atmospheric source inversion is to estimate the source term
$f(\x, t)$ on the right hand side of \eqref{eqnModelIntro} from indirect
measurements of the concentration $C(\x,t)$. We consider equation
\eqref{eqnParametric} with a set of given point sources at known
locations, but with unknown constant emission rates. Our aim is then to
estimate the emission rates $q_j$ from indirect measurements of
$C(\x,t)$. To this end, the source inversion problem is the reverse of
the forward problem \eqref{eqnParametric} and is thus an {\it inverse
  problem}.

We assume that measurements of concentration can be treated in terms of
the action of bounded linear operators on $C$. For example, a common
measurement in particulate dispersion studies is the total accumulated
deposition of pollutants in a region $R\subset\mathbb{R}^{2}$ of the
ground surface after some time $T$ has elapsed.  Such a measurement can
be written as
\begin{gather*}
  \int_{R} \int_{0}^{T} C(x,y,0,t)v_{\textit{set}}\,dt \, dx \, dy,
\end{gather*}
where $v_{\textit{set}}$ is the vertical settling velocity at the ground
surface for the pollutant of interest. The class of all possible linear
measurements is large and includes many commonly used methods in
practice ranging from short time average measurements such as Xact
ambient metal monitors \cite{cooper-environmental-2015} and Andersen
high-volume air samplers \cite{thermo-scientific-2015}, as well as
averaged long-time measurements accumulated in dustfall jars.

Suppose now that field measurements of $C$ are taken at $d$
locations in space and collected in a vector $\w \in
\mathbb{R}^d$. Since the measurements are linear, we can write
\begin{gather}
  \label{eqnLinear}
  \w = \mathcal{A}(\mytheta) \, \q,
\end{gather}
where $\mathcal{A}(\mytheta)$ is a $d\times n$ matrix that represents
the solution map of the PDE \eqref{eqnParametric} and depends nonlinearly
on the model parameters $\mytheta$. Source inversion requires finding
$\q$ given $\w$, and it is well known that this problem is ill-posed in
the sense of Hadamard \cite{enting1990inverse, isakov1990inverse}; that
is, small variations in the inputs ($\w$) result in large variations in
the solution ($\q$).  This is the case even if the true value of
$\mytheta$ is known exactly.

A variety of different approaches have been proposed in the literature
to solve source inversion problems~\cite{haupt-young-2008,
  rao-2007}. Lin and Chang \cite{lin2002relative} used an air trajectory
statistical approach to estimate the strength of different sources of
volatile organic compounds of anthropogenic origin. Stockie and Lushi
\cite{lushi2010inverse} used a Gaussian plume approximation to the
governing PDEs to estimate ground-level deposition of zinc from a
lead-zinc smelter, using linear least-squares to perform the source
inversion.  Skiba \cite{skiba2003method} solved the adjoint equation for
the advection-diffusion equation and used penalized least-squares to
invert the sources.

All of the methods just mentioned yield point estimates for the source
strengths but provide no direct measure of uncertainty in the estimated
parameters.  This drawback is overcome by the use of probabilistic
methods such as the Bayesian approach for solving inverse problems.
Such methods provide a robust setting for solving the inverse problem
and quantifying the uncertainties associated with the solution.  Sohn
et al.~\cite{sohn2002rapidly} developed an algorithm to obtain estimates
and uncertainties for the location and strength of pollutant sources in
buildings using data obtained from a COMIS simulation.  Hosseini and
Stockie \cite{hosseini2016bayesian} used a Gaussian plume model combined
with experimental data to estimate and quantify the uncertainty of
airborne fugitive emissions.  In a follow-up publication
\cite{hosseini2016airborne}, the same authors coupled a Bayesian
approach with a finite volume solver to estimate and quantify the
strength of airborne contaminants from a fixed number of point sources
at known locations.  Keats et al.~\cite{keats2007bayesian} obtained
probability distributions for source strengths and locations using
data from the standard ``Muck Urban Setting Test'' experiment.

The central contribution of this work is in proposing a method that
deals with the calibration parameters $\mytheta$ as unknowns. We infer
$\mytheta$ along with the vector of source strengths $\q$, which results
in a nonlinear inverse problem that is solved within the Bayesian
framework.  The problem of estimating the true value of the parameters
$\mytheta$ is referred to as {\it model calibration}.  Traditionally,
the process of tuning parameters in atmospheric dispersion models is
done empirically using Pasquill stability classes
\cite{turner1994workbook, seinfeld1998atmospheric} that are often chosen
heuristically.  Instead, we automatically estimate these parameters
using the information contained in the measured data $\w$.

\subsection{Outline}
The remainder of this article is organized as follows:
In
Section~\ref{sec:emulation} we introduce GP
emulators for approximation of computationally expensive models.
Section~\ref{sec:bayesianframework} is dedicated to the  Bayesian framework
for simultaneous calibration of the model and solution of the source
inversion problem. In Section~\ref{secIndustrialCase} we apply our  framework for
emulation and calibration  to an industrial case study involving
air-borne particulate emissions from a lead-zinc smelter in Trail,
British Columbia, Canada and compare our results with previous studies and engineering
estimates.

\section{Emulation of atmospheric dispersion models}
\label{sec:emulation}

In this section we introduce a framework for emulation of atmospheric dispersion models
using GPs.
Equation~\eqref{eqnParametric} captures the linear relationship between
source strengths $\q \in \mathbb{R}^n$ and concentration $C(\x,t)$.
Under the further assumption that the data $\w \in \mathbb{R}^d$ also
depends linearly on concentration $C(\x, t)$, we can write the map
$\q\mapsto \w$ in the linear form \eqref{eqnLinear}, although we note
that in general, $\mathcal{A}(\mytheta) \in \mathbb{R}^{d\times n}$ may
depend nonlinearly on $\mytheta \in \mathbb{R}^m$. Therefore, the map
$(\q, \mytheta) \mapsto \w$ can be split into two parts: for fixed
values of $\mytheta$ the map $(\cdot, \mytheta) \mapsto \w$ is linear,
whereas for a fixed $\q$ the mapping $(\q, \cdot) \mapsto \w$ is
nonlinear.  We will exploit this structure in the design of our source
inversion algorithm, using the fact that the linear map $(\cdot,
\mytheta) \mapsto \q$ can be dealt with efficiently, but that we will
need to approximate the nonlinear map $(\q, \cdot) \mapsto \w$.

To this end, we approximate the matrix
\begin{gather}
  \label{script-A-def}
  \mathcal{A}(\mytheta) = [\mathfrak{a}_{ij}(\mytheta)] \in \mathbb{R}^{d \times n},
\end{gather}
with another matrix 
\begin{gather*}
A(\mytheta) = [a_{ij}(\mytheta)] \in \mathbb{R}^{d \times n}, 
\end{gather*}
where each entry $\mathfrak{a}_{ij}: \mathbb{R}^m \mapsto \mathbb{R}$ of
$\mathcal{A}$ is approximated by the map $a_{ij}: \mathbb{R}^m \mapsto
\mathbb{R}$. For this purpose, we make use of GP
emulators of \cite{kennedy2001bayesian}.  Emulation with GPs is a
well-established method in dynamic computer experiments and machine
learning and we refer the reader to \cite{rasmussen2006gaussian,
  kennedy2001bayesian, o2006bayesian, conti2009gaussian} and the
references within for an introduction to this subject. Here, we outline
our method for approximating $\mathcal{A}$ and will not discuss the
theory of GP emulators with the exception of a few crucial definitions
and results.

\begin{definition}\label{dfnGP}
  A GP on $\mathbb{R}^m$ is a collection of
  real-valued random variables $\{g(\mathbf{x})\}_{\mathbf{x} \in
    \mathbb{R}^m}$, any finite collection of which have a joint Gaussian
  distribution.
\end{definition}
We denote a GP using the notation 
\begin{gather*}
  g \sim \mathcal{GP}(\bar{g}, \kappa), 
\end{gather*}
where the mean of the GP is
\begin{gather*}
  \bar{g}(\mytheta) = \mathbb{E}\: g(\mytheta) \qquad \forall \mytheta
  \in \mathbb{R}^m,  
\end{gather*}
and $\kappa$ is a positive semi-definite kernel with
\begin{gather*}
  \kappa( \mytheta, \mytheta') = \mathbb{E}\: (g(\mytheta) -
  \bar{g}(\mytheta))  (g(\mytheta') - \bar{g}(\mytheta')) \qquad
  \forall \mytheta, \mytheta' \in \mathbb{R}^m.    
\end{gather*}
We are primarily interested in isotropic kernels that satisfy
\begin{gather*}
  \kappa(\mytheta, \mytheta') = \kappa(|\mytheta - \mytheta'|),
\end{gather*}
from which it follows that if $\bar{g}$ is continuous and $\kappa$ is
continuous at zero, then the GP is also mean square continuous
\cite{rasmussen2006gaussian}.

Let $\{\mytheta_k\}_{k=1}^K$ be a collection of {\it design points} in
the parameter space $\mathbb{R}^m$ for some fixed $K > 0$, and let
$\{\mathbf{e}_j\}_{j=1}^n$ be the unit coordinate basis vectors in
$\mathbb{R}^n$. Then for each index pair $(i,j)$ define a GP
\begin{gather*}
  g_{ij} \sim \mathcal{GP}(\bar{g}_{ij}, \kappa_{ij}),
\end{gather*}
subject to the constraint  
\begin{gather*}
  g_{ij}(\mytheta_k) = \mathfrak{a}_{ij}(\mytheta_k) =
  (\mathcal{A}(\mytheta_k)\mathbf{e}_j)_i 
  \qquad \text{for } k = 1,\ldots, K.
\end{gather*}
That is, each $g_{ij}$ interpolates the $\mathfrak{a}_{ij}$ at points
$\{\mytheta_k\}_{k=1}^K$.  We can now identify $g_{ij}$ using well-known
identities for conditioning of Gaussian random variables
\cite[Sec.~2.2]{rasmussen2006gaussian}.  Let $\Theta := (\mytheta_1,
\ldots, \mytheta_K)^T$ and define the following matrices and vectors
\begin{equation*}
  \begin{aligned}
    & G_{ij}( \Theta, \Theta) \in \mathbb{R}^{K\times K}, 
    \qquad &&[G_{ij}(\Theta, \Theta)]_{k\ell} := 
    \kappa_{ij}(\mytheta_k, \mytheta_\ell), 
    \\
    & G_{ij}(\mytheta, \Theta) \in \mathbb{R}^{1 \times K}, 
    \qquad &&[G_{ij}(\mytheta, \Theta)]_{\ell} := 
    \kappa_{ij}(\mytheta, \mytheta_\ell), 
    \\
    & G_{ij}(\Theta, \mytheta) \in \mathbb{R}^{K \times 1}, 
    \qquad &&[G_{ij}(\Theta, \mytheta)]_{k} := 
    \kappa_{ij}(\mytheta_k, \mytheta), 
    \\ 
    & \mathbf{g}_{ij} \in \mathbb{R}^{K \times 1}, 
    \qquad && [\mathbf{g}_{ij}]_k = \mathfrak{a}_{ij}(\mytheta_k), 
  \end{aligned}
\end{equation*}
for $k,\ell = 1, \ldots, K$.  Then we have
\begin{gather*}
  g_{ij}(\mytheta) \sim \mathcal{GP}(\bar{g}_{ij}(\mytheta),
  \sigma_{ij}(\mytheta)) , 
\end{gather*}
where 
\begin{equation*}
  \begin{aligned}
    \bar{g}_{ij}(\mytheta) & =  G_{ij}(\mytheta, \Theta) 
    G_{ij}(\Theta, \Theta)^{-1} \mathbf{g}_{ij}, \\ 
    \sigma_{ij}(\mytheta)  & = \kappa_{ij}
    (\mytheta, \mytheta) - G_{ij}(\mytheta, \Theta) \, 
    \left[G_{ij}(\Theta, \Theta)^{-1}\right] G_{ij}(\Theta, \mytheta),
  \end{aligned}
\end{equation*}
for all $\mytheta \in \mathbb{R}^m$. Since the mean $\bar{g}_{ij}$
interpolates the data $\mathfrak{a}_{ij}(\mytheta_k)$, it is also a good
candidate for approximating $\mathfrak{a}_{ij}$, and so we take
\begin{gather*}
  a_{ij}(\mytheta) = \bar{g}_{ij}(\mytheta) 
  \qquad \forall \mytheta \in \mathbb{R}^m. 
\end{gather*}
The advantage of using this approach for interpolation is that it is
possible to assess the uncertainty and quality of the emulator via the
covariance operator $\sigma(\mytheta)$.

The variance at a point $\mytheta_{i}$ is given by the term
$\sigma_{ii}(\mytheta)$ and is strongly influenced by the spatial
distribution of the points $\{\mytheta_{k}\}_{k=1}^{K}$
\cite{johnson1990minimax}.  Thus, in order to emulate
$\mathcal{A}(\mytheta)$, it is crucial to choose the points
$\mytheta_{k}$ in such a way that the uncertainty of the emulator is
minimized.  A popular method for choosing the design points is known as
a space-filling design \cite{johnson1990minimax}, instances of which
include the maximum entropy design \cite{maxEntDesign}, uniform design
\cite{Fang}, and Latin Hypercube design (LHD) \cite{LHD}. Following
Jones et al.~\cite{jones} we use a combination of the LHD and maximin
designs, which is introduced in \cite{johnson1990minimax} and shown to
outperform most space-filling design methods for GP
interpolation.

We now briefly outline our space-filling design procedure.  The main
idea in the LHD is to distribute design points such that low-dimensional
projections of the points do not overlap. One issue with LHD is that it
does not have desirable space-filling properties and so leads to large
uncertainties in the emulator in high dimensions. Thus, LHD is often
used as an initial condition for other space-filling design
methodologies and so here we complement the LHD with an approximate
maximin design \cite{johnson1990minimax}.
 
The idea behind the maximin design is as follows.  Let
$T\subset\mathbb{R}^{n}$ be the subset of parameter space in which we
wish to construct our design points and consider all subsets $S$ of $T$
with finite (fixed) cardinality, say $|S|=k$.  A maximin design $S^{o}$
is a set of points that satisfies
\begin{gather}
  \label{eqnmaximin}
  \max_{S\subset T,\text{ }|S|=k} \; \min_{s,s'\in S}d(s,s') =
  \min_{s,s'\in S^{o}}d(s,s'), 
\end{gather}
where $d$ is the Euclidean metric. The optimization problem in equation
\eqref{eqnmaximin} is not easily solvable, which is why we resort to
metaheuristic optimization procedures.  In particular, we use the
particle swarm algorithm \cite{arora2015optimization}, initiated with an
LHD. We then apply a number of iterations of particle swarm that is
chosen depending on our available computational budget. In
Figure~\ref{figSpaceFilling} we show an example of the final design
after 10,000 steps of a particle swarm algorithm.

\section{The Bayesian inversion framework}
\label{sec:bayesianframework}

We introduce the Bayesian framework for calibration and solution of
inverse problems in this section.
 The Bayesian approach
combines the data $\w$ together with prior knowledge to give a posterior
probability distribution on the source strengths $\q$ and 
calibration parameters $\mytheta$.  Following Bayes' rule, the posterior
distribution of $(\mytheta, \q)$ is given by
\begin{gather}
  \label{eqnBayesChp4}
  \post(\mytheta,\q|\w) = \frac{\like(\w|\mytheta,\q)
    \, \prior(\mytheta,\q)}{Z(\w)}, 
\end{gather}
where $\like(\w|\mytheta,\q)$ is the \textit{likelihood probability} of
$\w$ given both $\mytheta$ and $\q$, which represents the probability of
observing $\w$ given fixed values $\q$ and $\mytheta$.  The probability
distribution $\prior(\mytheta,\q)$ is called the \textit{prior
  probability} that expresses prior knowledge of the possible values
of $\q$ and $\mytheta$, before observing any experimental data.  Finally,
$Z(\w)$ is a constant that normalizes the posterior probability so that
it integrates to one, leading to the requirement that
\begin{gather}
  \label{eqnPropConst}
  Z(\w) = \int \like(\w|\mytheta,\q) \, \prior(\mytheta,\q)\,d\mytheta\,
  d\q. 
\end{gather}

We next assume that the linearity requirement \eqref{eqnLinear} holds
for physical measurements up to some additive Gaussian measurement
noise, so that
\begin{gather}
  \label{eqnNoisy}
  \w = \mathcal{A}(\mytheta) \q + \epsilon \quad \text{with} \quad 
  \epsilon \sim \mathcal{N}(0, \Sigma). 
\end{gather}
Here, $\epsilon$ is a normally distributed random vector with mean zero
and covariance $\Sigma \in \mathbb{R}^{d\times d}$, which is assumed
positive-definite. The measurement noise $\epsilon$ models the
deviations in the measurements $\w$ due to missing physics and
uncertainty in measurements. Under this assumption, it can be readily
shown that \cite[Sec.~3.2.1]{Somersalo}
\begin{gather}
  \label{eqnLike}
  \like(\w|\mytheta,\q) = \frac{1}{(2\pi \text{det}
    \Sigma)^{\frac{1}{2}}}\exp\left(-\frac{1}{2} \left\|
    \Sigma^{-1/2}(\mathcal{A}(\mytheta)\q-\w ) \right\|^{2} \right). 
\end{gather}

The choice of the prior distribution depends on previous knowledge about
the model parameters and so is problem specific. For atmospheric
dispersion models, the parameter $\mytheta$ often depends on atmospheric
conditions, whereas the emission rates in the vector $\q$ depend on the
physical processes that generated the emissions.  Consequently, it is
reasonable to assume a priori that $\q$ and $\mytheta$ are statistically
independent, so that
\begin{gather}
  \label{eqnPrior}
  \prior(\mytheta,\q) = \prior(\mytheta) \, \prior(\q).
\end{gather}
In most cases, one chooses a prior for $\q$ that imposes a positivity
constraint on the emission rates, whereas priors on $\mytheta$ typically
reflect acceptable ranges of calibration parameters.  Beyond these
assumptions, the choice of the prior density for $\mytheta$ and $\q$
must be made on a case by case basis.  In
Section~\ref{secIndustrialCase}, we present a particular prior
distribution for $(\mytheta, \q)$ in the context of an industrial case
study.

A key step in our source inversion approach is approximating the
posterior probability $\post$ through an approximation for the map
$\mathcal{A}$. Let $A$ be the GP emulator for $\mathcal{A}$ as outlined
in Section~\ref{sec:emulation}, then we can approximate the likelihood
$\like$ in \eqref{eqnLike} using
\begin{gather}
  \label{approx-eqnLike}
  \widehat{\like}(\w|\mytheta,\q)=\frac{1}{(2\pi \text{det}
    \Sigma)^{\frac{1}{2}}}\exp\left(-\frac{1}{2} \left\|
    \Sigma^{-1/2}({A}(\mytheta)\q-\w ) \right\|^{2}\right). 
\end{gather}
Substituting this expression into Bayes' rule \eqref{eqnBayesChp4}
yields the approximate posterior
\begin{gather}
  \label{approx-eqnBayesChp4}
  \widehat{\post}(\mytheta,\q|\w) =
  \frac{\widehat{\like}(\w|\mytheta,\q) \,
    \prior(\mytheta,\q)}{\widehat{Z}(\w)},  
\end{gather} 
where the constant $\widehat{Z}$ is defined similar to
\eqref{eqnPropConst}.  If the map $A$ approximates $\mathcal{A}$
closely, then we expect $\widehat{\post}$ to be close to $\post$. We
refer the reader to \cite{stuart2018posterior} and references therein
for a detailed analysis of the errors introduced through GP emulation of
forward maps in Bayesian inverse problems.

With $\widehat{\post}$ in hand we may now compute point value estimators
for the parameters $\mytheta$ and source strengths $\q$.  Common choices
of point estimates are \cite{Somersalo}:
\begin{alignat}{4}
  \label{eqnpointestimates}
  \text{Maximum a posteriori:} &\qquad &
  (\mytheta,\q)_{\textit{MAP}} &=
  \argmax_{\mytheta,\q} \, \widehat{\post}(\w|\mytheta,\q), 
  \\ 
  \text{Conditional mean:} &&
  (\mytheta,\q)_{\textit{CM}} &=
  \int(\mytheta,\q) \, \widehat{\post}(\mytheta,\q|\w) \,d\mytheta\, d\q,  
  \\
  \text{Maximum likelihood:} &&
  (\mytheta,\q)_{\textit{ML}} &=
  \argmax_{\mytheta,\q} \, \widehat{\post}(\mytheta,\q|\w).
\end{alignat}
The uncertainty in the choice of point estimate
$(\mytheta^{\ast},\q^{\ast})$ can be assessed by computing the
covariance matrix
\begin{gather*}
  \int \left( (\mytheta,\q)-(\mytheta^{\ast},\q^{\ast}) \right) \otimes
  \left( (\mytheta,\q)-(\mytheta^{\ast},\q^{\ast}) \right) \,
  \widehat{\post}(\mytheta,\q|\w) \, d\mytheta \,d\q.
\end{gather*}
The above estimators may of course also be computed using the true
posterior ${\post}$, but if $\widehat{\post}$ is close to $\post$ in an
appropriate sense (such as the total variation metric) then one expects
the point estimators under $\widehat{\post}$ to approximate their
counterparts under $\post$ \cite{stuart2018posterior}.

In general, the integrals involved in computing point estimates or
covariances are not analytically tractable and so it is necessary to
resort to numerical methods for their estimation. Since these are
high-dimensional integrals, quadrature-based approaches are unsuitable,
thus it is necessary to use Markov chain Monte Carlo (MCMC) integration
techniques. In this paper we use the adaptive MCMC algorithm of
\cite{haario2001}, which is outlined below in Algorithm~\ref{algMH}.
Specific values of the algorithmic parameters ($\beta_1$, $\gamma_1$,
$\gamma_2$, $\gamma_3$) that are tailored to our case study are provided
later in Section~\ref{sec:source-inversion}.
 
\begin{algorithm}
  \caption{Adaptive Metropolis-Hastings Algorithm.}\label{algMH}
  \begin{algorithmic}[1]
     \State Choose $(\mytheta_1, \q_1) \in \mathbb{R}^d$ in the support of
    $\post(\mytheta,\q|\w)$ and fixed parameters $\beta \in (0,1)$,
    $\gamma_1, \gamma_2, \gamma_3 \in (0,\infty)$, $N\in \mathbb{N}$. 
    \For{$j=2:N$}
    \If{$j\leq 2d$}
    \State Draw $u$ from $\mathcal{N} \left(
      (\mytheta_{j},\q_{j}), \frac{\gamma_1}{d} I_{d\times d} \right)$. 
    \Else
    \State Estimate the empirical covariance matrix $\Sigma_{j}$ using $\{ (\mytheta_k, \q_k) \}_{k=1}^j$. 
    \State Draw $u$ from $(1-\beta) \, \mathcal{N} \left(
      (\mytheta_{j},\q_{j}), \frac{\gamma_2}{d}\Sigma_{j} \right) + 
    \beta\mathcal{N}\left( (\mytheta_{j},\q_{j}),
      \frac{\gamma_3}{d}I_{d\times d} \right)$. 
    \EndIf
    \State Propose $(\tilde{\mytheta}_j,\tilde{\q}_{j}) = 
    (\mytheta_{j-1},\q_{j-1}) + u$.
    \State Compute $\delta = \min
    \left(1,\frac{\widehat{\post}(\mytheta_{j},\q_{j}|\w)}
      {\widehat{\post}(\mytheta_{j-1},\q_{j-1}|\w)} \right)$. 
    \State Draw $w$ from $U([0,1])$.
    \If{$w<\delta$}
    \State $(\mytheta_{j},\q_{j}) =
    (\tilde{\mytheta}_{j},\tilde{\q}_{j}) \qquad\qquad$ (Accept the move) 
    \Else
    \State $(\mytheta_{j},\q_{j})=(\mytheta_{j-1},\q_{j-1})\qquad$ 
    (Reject the move) 
    \EndIf
    \EndFor
  \end{algorithmic}
\end{algorithm}



\section{An industrial case study}
\label{secIndustrialCase}

\hl{ In this section we apply our Bayesian framework to the study of
  dispersion of airborne zinc particles from four point sources located
  within the area surrounding a lead-zinc smelter in Trail, British
  Columbia, Canada.  In Section~\ref{sec:sensitivity-analysis} we
  perform a sensitivity analysis that reveals which parameters are most
  informed by the data. Section~\ref{sec:construct-emulator} is
  dedicated to validating the emulator while
  Section~\ref{sec:source-inversion} contains the main results for the
  industrial case study including comparisons with previous studies.
  Finally, Sections~\ref{sec:effect-prior-choice} and
  \ref{sec:effect-emul-qual} provide further validation of our source
  inversion results and study the dependence of emission estimates
  and uncertainties on the choice of the priors and quality of the
  emulator respectively. The methodology of this section was validated
  in \cite[Sec.~2]{garcia2018parameter} for a synthetic example, and the
  dependence of the posterior on the number of measurements and choice
  of the prior distribution was studied in detail.  }

\begin{figure}[htbp]
  \centering
  \includegraphics[width=0.6\textwidth]{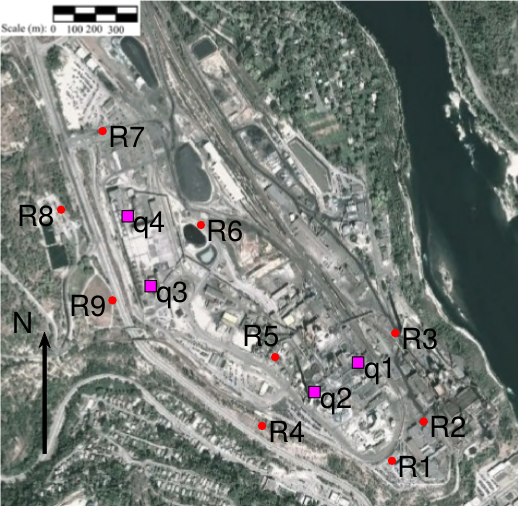}
  \caption{Aerial photograph of the lead-zinc smelter in Trail, BC,
    Canada. The square points labelled 'qn' represent the locations of
    the four zinc sources, while the circular points 'Ri' denote the
    nine measurement devices.}
  \label{figAerial}
\end{figure}

Our ultimate goal is to estimate the
contribution that each source makes to the total zinc released into the
atmosphere by the smelting operation. We begin by outlining the details of our
model and the parameters to be calibrated.
We have access to monthly
cumulative measurements of zinc particulate depositions at nine separate
locations as well as horizontal wind field velocity data at a
meteorological station located near the sources. We also have access to
engineering estimates of the yearly-averaged emission rates obtained
from independent engineering studies based on process control arguments.
An aerial photograph of the industrial site, showing the locations of
all sources and measurement devices, is provided in
Figure~\ref{figAerial}. The sources are labelled $q_{1}$ to $q_{4}$, and
deposition measurements are designated $R_{1}$ to $R_{9}$.  This same
emission scenario has been studied using a linear least squares approach
based on a Gaussian plume approximation~\cite{lushi2010inverse}, and
also by performing a finite volume discretization of the governing
equation~\eqref{eqnParametric}~\cite{hosseini2016airborne}.

We apply the atmospheric dispersion model \eqref{eqnParametric},
as adapted in~\cite{hosseini2016airborne} for this specific zinc
smelter site with $n=4$ sources:
\begin{gather}
  \label{eqnParametricFinal}
  \frac{\partial C(\x,t)}{\partial t} +
  \textbf{L}(p,z_{0},z_{i},L,z_{\textit{cut}}) C(\x,t) 
  = \sum_{j=1}^{4} q_{j}\delta(\x-\x_{j}). 
\end{gather}
The differential operator $\textbf{L}$ depends non-linearly on five
parameters as we shall explain next.  The parameter $p$ comes from the
assumption of a power-law distribution in the vertical component of the
wind profile $\vv$; that is, if $z$ denotes the height above ground
level $(x,y,0)$, then the two horizontal velocity components $v_{x}$ and
$v_{y}$ are assumed to be function of $z$ and $t$ only and satisfy
\begin{gather*}
\|(v_{x}(z,t),v_{y}(z,t))\|_{2}=v_{r}(t)\left(\frac{z}{z_{r}}\right)^{p},
\end{gather*}
where $v_{r}(t)$ is the wind speed at some reference height $z_{r}$.
\hl{Note that the Euclidean norm appears on the left hand side  since the power law describes
  the magnitude of the horizontal wind field.}
The
exponent $p$ depends on factors such as surface roughness and
atmospheric stability class. For more details about the power law model
for the wind velocity the reader is referred to
\cite[Sec.~16.4.5]{seinfeld1998atmospheric}.  The other four parameters are related
to the entries in the diagonal diffusivity matrix
$\textbf{D}=\operatorname{diag}(D_{11}, D_{22}, D_{33})$ in
equation~\eqref{eqnParametric}.  Following
\cite{seinfeld1998atmospheric}, the vertical diffusion coefficient
satisfies to a good approximation
\begin{gather}
  \label{eqnEddyVertical}
  D_{33}=\frac{\kappa v_{\ast} z}{\phi(z/L)},
\end{gather}
where $\kappa$ is the \textit{von K{\'a}rm{\'a}n} constant whose value can be
set to $0.4$ in practical scenarios. The function $\phi$ in the
denominator is taken to be a piecewise continuous function
\begin{gather*}
  \phi\left(\frac{z}{L}\right)=\left\{
    \begin{array}{ll}				
      1+4.7\frac{z}{L},                 &\mbox{if }\frac{z}{L}\geq 0, \\
      (1-15\frac{z}{L})^{-\frac{1}{2}}, &\mbox{if }\frac{z}{L}<0,
    \end{array}
  \right.
\end{gather*}
where $L$ is called the \textit{Monin-Obukhov length}. The parameter
$v_{\ast}$ is known as the \textit{friction velocity} and is represented by 
\begin{gather*}
  v_{\ast}(t) = \frac{\kappa v_{r}(t)}{\ln(\frac{z_{r}}{z_{0}})}, 
\end{gather*}
where $v_{r}(t)$ is the wind velocity at some reference height. The
variable $z_{0}$ is called the \textit{roughness length} and depends on
both terrain and surface type.  The remaining horizontal diffusion
coefficients satisfy $D_{11}=D_{22}$ and are assumed to be independent
of height $z$. A correlation that is commonly employed for these
parameters is \cite{seinfeld1998atmospheric}
\begin{gather*}
  D_{11}=D_{22}\approx
  \frac{v_{\ast}z_{i}^{\frac{3}{4}}(-\kappa L)^{-\frac{1}{3}}}{10},
\end{gather*}
where the parameter $z_{i}$ is called the \textit{mixing layer height}.

The governing PDE~\eqref{eqnParametricFinal} is complemented with 
boundary conditions defined over the space-time domain
$\mathbb{R}^{2}\times[0,\infty)\times(0,T)$ for some fixed
$T\in\mathbb{R}$.  For this purpose, we impose the following conditions
at infinity
\begin{gather*}
  c(\x,t)\rightarrow 0\qquad\text{as }\|\x\|\rightarrow\infty,
\end{gather*}
along with a flux (or Robin) boundary condition at the ground surface
\hl{
\begin{gather}
  \label{eqnRobinBoundary}
  \left.\left(v_{\textit{set}}c - D_{33}\frac{\partial
        c}{\partial z}\right)\right\rvert_{z=0} =
  \left.v_{\textit{dep}}c\right\rvert_{z=0},  
\end{gather} 
where $v_{\textit{set}}=0.0027\;m/s$ is the settling velocity of zinc
particles computed using Stokes' law \cite[Sec.~9.3]{seinfeld1998atmospheric} and $v_{\textit{dep}}=0.005\; m/s$ is the corresponding
deposition velocity.} Observe that at the bottom boundary $z=0$, the
vertical diffusivity satisfies $D_{33}(0)=0$ from
\eqref{eqnEddyVertical} which leads to an inconsistency in the Robin
boundary condition; therefore, we define a cut-off length
$z_{\textit{cut}}\gtrapprox 0$ below which $D_{33}$ is set to the
non-zero constant value $D_{33}(z_{\textit{cut}})$.

Equation~\eqref{eqnParametricFinal} is solved numerically using the
finite volume algorithm described \cite{hosseini2016airborne}.  The code
takes as input values of the wind speed $v_r(t)$, parameters
$p,z_{0},z_{i},L,z_{\textit{cut}}$, and source strengths $q_{j}$, and
returns as output the zinc concentration at any point $(x,y,z)$
and time $0<t<T$. However, as mentioned earlier, we do not have
experimental measurements of concentration but rather total depositions
at the nine sites $R_{1},\ldots,R_{9}$ (see Figure~\ref{figAerial}). We
model the depositions using the integral
\begin{gather}
  \label{eqnw}
  w(x,y,T) = \int_{0}^{T} c(x,y,0,t) \, v_{{\textit set}} \,dt.
\end{gather}
  Since the nine
deposition measurements $R_{i}$ were obtained from roughly cylindrical 
dustfall jar collectors, we can readily approximate the total deposition
over the time interval $(0,T]$ at site $i=1,\ldots,9$ using 
\begin{gather}
  \label{eqnDepositionIntegral}
  w_i = \int_{\textit{jar}} w(x,y,T) \,dx \,dy 
  \approx w(x_{i},y_{i},T) \Delta A,
\end{gather}
where $\Delta A\approx 0.0206\;m^2$ is the cross-sectional area of each
dustfall jar. We then collect the zinc deposition measurements into a
vector $\w=(w_i)^T$ and from now on, we assume that the time interval
$T$ over which particulates are allowed to accumulate in the dustfall
jars is one month.

With the industrial case fully specified our source inversion framework
can be implemented.  We make the following observations:
\begin{itemize}
\item The vector $\w$ is composed of $d=9$ measurements, one from each
  jar, so that $\w\in\mathbb{R}^{9}$.
\item The number of model parameters is $m=5$, which are represented by 
  the vector $\mytheta = (p,z_{0},z_{i},L,z_{\textit{cut}})^T$. 
\item There are $n=4$ sources stored in the vector
  $\q=(q_{1},q_{2},q_{3},q_{4})^T$.
\end{itemize}
Each measurement in the vector $\w$ may then be viewed as a map
\begin{gather*}
  w_{i}:\mathbb{R}^{5} \times
  \mathbb{R}^4\rightarrow [0,\infty). 
\end{gather*}
In other words, the number $w_{i}(\mytheta^{\ast},\q^{\ast})$ is the
deposition at site $R_{i}$ predicted by the finite volume solver
described in \cite{hosseini2016airborne}, given a particular choice of
parameters $\mytheta^{\ast}$ and source strengths $\q^{\ast}$.  In the
framework of Section~\ref{sec:emulation} (namely, equations
\eqref{eqnLinear} and \eqref{script-A-def}) we write
\begin{gather}
  \label{eqnyimaps}
  w_{i}(\mytheta^{\ast},\q^{\ast}) = \sum_{j=1}^{4}
  \mathfrak{a}_{ij}(\mytheta^{\ast}) \, q^{\ast}_{j},
\end{gather}
where each of the nonlinear functions
$\mathfrak{a}_{ij}(\mytheta^{\ast})$ can be approximated using a GP
emulator $a_{ij}(\mytheta^{\ast})$.

\subsection{Sensitivity analysis}
\label{sec:sensitivity-analysis}

Before solving the source inversion problem, we first perform a sensitivity
analysis in order to check whether the dimensionality of the parameter
space can be reduced. The goal of a sensitivity analysis is to assess
the relative importance of each of the variables in a model
\cite{saltelli2000sensitivity}.  In this paper, the sensitivity measure
we employ is known as the Sobol total index \cite{sobol1993sensitivity,
  sullivan2015introduction, saltelli2000sensitivity} which captures the
variance of the model outputs due to changes in each of its inputs and
is a measure of the relative importance of each input. More
specifically, the Sobol total index ranges from a lower limit of $0$,
denoting a variable that has no effect on the output, to an upper limit
of $1$, indicating that all variability in the model is explained by the
variable in hand.

We perform our sensitivity analysis using the R package
\pkgname{sensitivity} \cite{sensitivity}, which uses GPs
as surrogates for the function of interest. In order to account for the
variability introduced by using different isotropic kernels (see
Definition~\ref{dfnGP}), we calculate Sobol indices using the following
choices of the kernel $\kappa$:
\begin{itemize}
\item Exponential:
  $\kappa(|\mytheta-\mytheta'|) = r_{1}\exp\left(
    -\frac{|\mytheta'-\mytheta|}{r_{2}} \right)$,\\ 
\item Squared exponential: $ \kappa(|\mytheta-\mytheta'|) =
  r_{1}\exp \left( -\frac{|\mytheta-\mytheta'|^2}{2r_{2}} \right), $
\item Mat{\'e}rn $\frac{3}{2}$: $\kappa(|\mytheta-\mytheta'|) = r_{1} \left(1
    + \frac{\sqrt{3}|\mytheta-\mytheta'|}{r_{2}}\right) \exp\left(
    -\frac{\sqrt{3}|\mytheta-\mytheta'|}{r_{2}} \right)$,
\item Mat{\'e}rn $\frac{5}{2}$: $\kappa(|\mytheta-\mytheta'|) = r_{1} \left(1
    + \frac{\sqrt{5}|\mytheta-\mytheta'|}{r_{2}} + \frac{5}{3}
  \left(\frac{|\mytheta-\mytheta'|}{r_{2}}\right)^{2}\right) \exp\left( 
    -\frac{\sqrt{5}|\mytheta-\mytheta'|}{r_{2}} \right)$.
\end{itemize}
In the above,  $r_{1},r_{2}\in\mathbb{R}$ are tuning parameters.
These kernels are commonly used in 
GP regression and are readily implemented in the R packages
\pkgname{sensitivity} and \pkgname{DiceKriging}. For 
other possible choices of kernels see \cite[Sec.~4]{rasmussen2006gaussian}

\begin{figure}[htbp]
  \centering
  \includegraphics[width=0.9\textwidth]{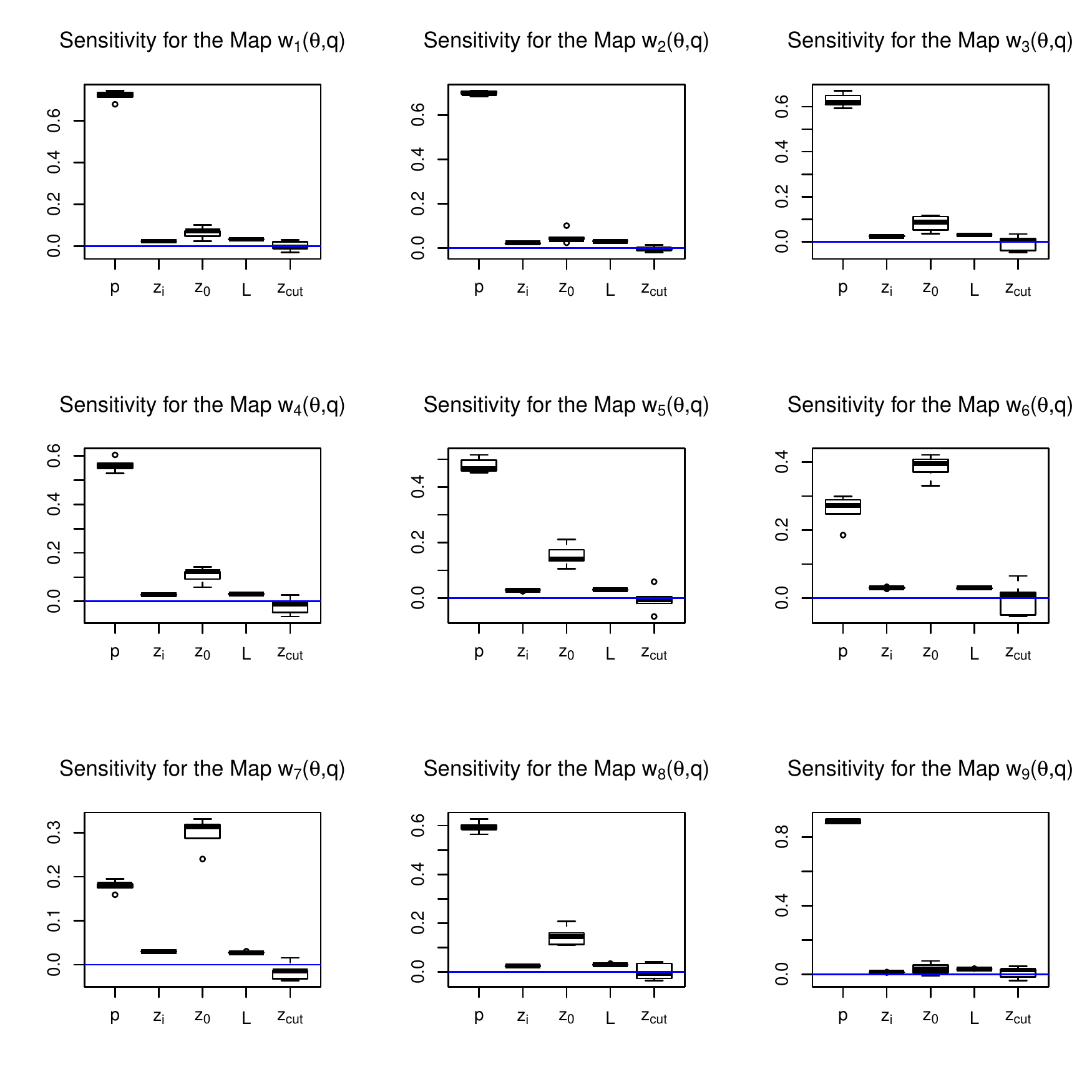}
  \caption{The Sobol total indices of the model parameters $\mytheta$
    for the maps $w_{i}$ defined in \eqref{eqnyimaps}. Errors bars
    indicate variation between estimated value of Sobol total indices
    using different kernels.}
  \label{figPaperSensitivity}
\end{figure}

The sensitivity indices computed using the four different kernels are
summarized in the boxplots shown in Figure~\ref{figPaperSensitivity}.
The size of each box represents the interquartile range (IQR).  \hl{The
upper bar attached to each box is at the value
\begin{gather*}
  \min \left\{ \max(S),E_{3}+\frac{3}{2}IQR \right\},
\end{gather*}
and the lower bar is located at 
\begin{gather*}
  \max \left\{\min(S),E_{1}-\frac{3}{2}IQR\right\},
\end{gather*}
}
where $S$ is the set of the four indices (corresponding to the four
kernels) calculated for each $w_{i}$ for $i\in\{1,\ldots,9\}$, and
$E_{1}$ and $E_{3}$ are the first and third quartiles of $S$
respectively~\cite{boxplot}.

From these results it is clear that the variability in output of each
map $w_{i}$ is captured in large part by the two parameters $p$ and
$z_0$.  Although the sensitivity index for $L$ is relatively small, we
still retain this parameter because it is known to be closely related to
$z_{0}$~\cite[Chap. 19]{seinfeld1998atmospheric}.  The remaining
``unimportant'' parameters, $z_i$ and $z_{\textit{cut}}$, may then be
fixed at suitable constant values (we used $z_i = 100$ and
$z_{\textit{cut}} = 2$, following \cite{hosseini2016airborne}).
This reduces the dimensionality of the parameter space from $5$ to $3$ and
allows us to redefine our parameter vector as
\begin{gather*}
  \mytheta := (p,z_{0},L).
\end{gather*}
\hl{We highlight that the negative values of the total Sobol index for
  $z_{\textit{cut}}$ appearing in Figure~\ref{figPaperSensitivity} are
  not physical but rather arise
  due to numerical errors in the Monte Carlo integration method used in
  the \pkgname{sensitivity} package.}

\begin{figure}[htbp]
  \centering
  \includegraphics[width=0.7\textwidth]{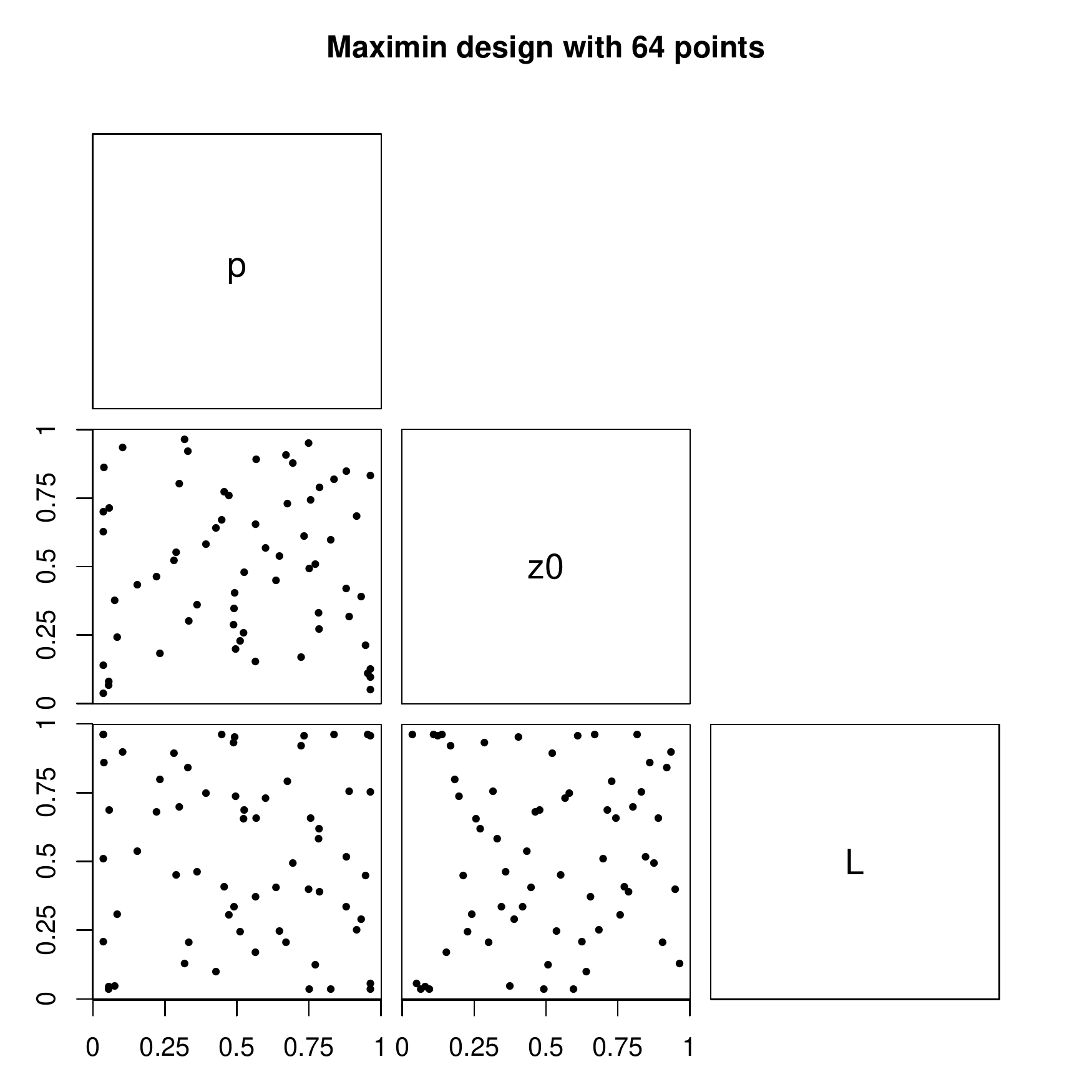}
  \caption{Pairwise two-dimensional projections of the space-filling
    design with 64 points for emulating $\mathcal{A}$.  Results were
    obtained by initializing the particle swarm algorithm using a LHD
    and iterating for 10,000 steps. The range of each parameter is
    normalized between $0$ and $1$ for easy visualization.  The scatter
    plot in the second row depicts $z_0$ vs.\ $p$. The plots on the
    bottom row show $L$ vs.\ $p$ (left) and $L$ vs.\ $z_0$ (middle).}
  \label{figSpaceFilling}
\end{figure}

We next identify the support of the prior distribution, which consists
simply of the interval of allowable values for each unknown parameter.
Every emission rate must be a non-negative number, so that 
\begin{gather*}
  q_{j}\in [0,\infty),\qquad\text{for } j=1,2,3,4.
\end{gather*}
For the remaining parameters $\pars$, we use the physically-motivated
ranges listed in Table~\ref{tabNewRanges}, which are slightly larger
than the suggested ranges typically cited in the literature (see e.g.,
\cite[Ch.~19]{seinfeld1998atmospheric}) to allow for more flexibility.
Consequently, each map $w_{i}(\mytheta,\q)$ has a domain
\begin{gather*}
  \text{dom}(w_{i}) = [0,0.6]\times[0,3]\times[-600,0] \times
  [0,\infty)^{4}. 
\end{gather*}
With this information in hand, we are now in a position to create a
space-filling design needed to construct the GP emulator as explained in
Section~\ref{sec:emulation}.  Due to computational budget
constraints, we choose 64 points to construct our emulators. We display
 2D projections of the corresponding maximin space-filling design for
the three relevant parameters in Figure~\ref{figSpaceFilling}. 
\begin{table}[htbp]
  \centering
  \caption{Physically reasonable ranges for the parameters $\pars$.}
  \label{tabNewRanges}
  \begin{tabular}{|c|c|c|}
    \hline 
    Parameter (units) & Symbol & Range\tabularnewline
    \hline 
    Velocity exponent  & $p$ & ${[}0,0.6{]}$\tabularnewline
    \hline 
    Roughness length (m) & $z_{0}$ & ${[}0 ,3{]}$\tabularnewline
    \hline 
    Monin-Obukhov length (m) & $L$ & ${[}-600,0{]}$\tabularnewline
    \hline 
  \end{tabular}
\end{table}

\subsection{Construction and validation of the emulator}
\label{sec:construct-emulator}

We now proceed as outlined in Section~\ref{sec:emulation} to construct
the emulator for maps $\mathfrak{a}_{ij}$ using GPs. The 64 different combinations of parameters $\pars$ obtained from
the space-filling design in the previous section are used to train the
GP. To create the GPs for each of the 36 maps
$\mathfrak{a}_{ij}(p,z_{0},L)$, we used the R package
\pkgname{DiceKriging} \cite{dupuy2015dicedesign}. Unlike the sensitivity
analysis in the previous sections where four different kernels were
computed, the GPs here are constructed using only the isotropic squared
exponential kernel.
 
\begin{figure}[htbp]
  \centering
  \includegraphics[width=0.5\textwidth]{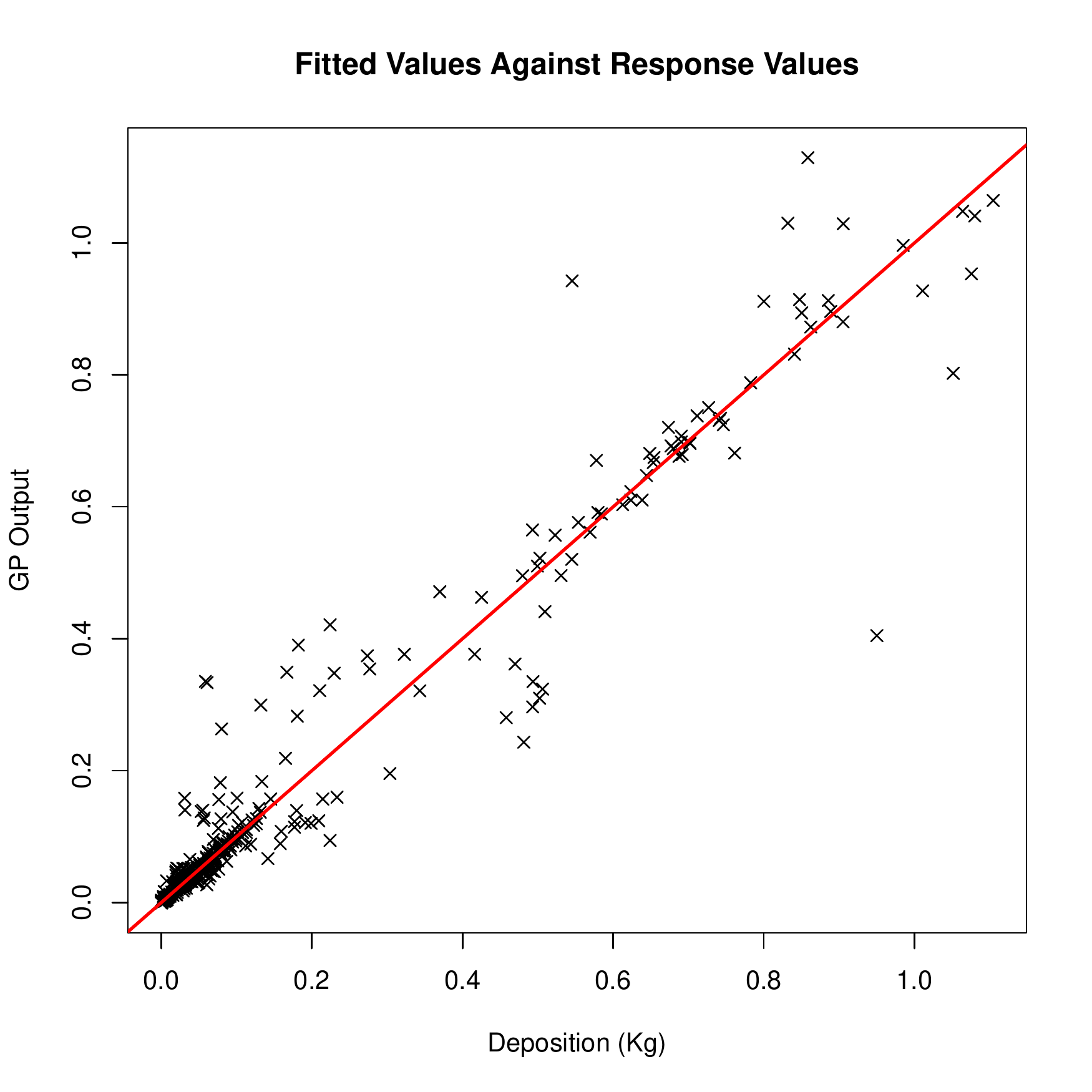}
  \caption{Comparison between deposition values predicted by the finite
    volume solver and the GP emulator using LOOCV for each of the
    $64\times 9$ data points obtained from the finite volume solver.}
  \label{figGPTrained}
\end{figure}

The \pkgname{DiceKriging} package chooses a maximum likelihood estimate
for the tuning parameters $r_1, r_2$ for each map
$\mathfrak{a}_{ij}$. To assess the quality of each emulator we perform a
``leave one out'' cross validation (LOOCV)
\cite[Sec~1.4]{murphy2012machine} on each of the nine sites
$R_{1},R_{2},\ldots,R_{9}$ for the 64 design points shown in
Figure~\ref{figSpaceFilling}; that is, we run the finite volume solver
64 times and save the predicted deposition value for each sites,
yielding a total of $64\times 9$ points or predictions from the finite
volume solver.  For each point ``left out'' in the LOOCV, we plot in
Figure~\ref{figGPTrained} the deposition obtained from the finite volume
solver versus the predicted value from the corresponding trained GP.
This figure shows that the output of the GP is closely matched to the
output of the finite volume solver at most design points, with the
exception of a small number of outliers.

\subsection{Source inversion}
\label{sec:source-inversion}

The GP emulator may now be employed to estimate the emission rates
$q_1,\ldots,q_4$, for the actual physical dataset which consists of
dustfall jar measurements at the nine locations marked
$R_{1},\ldots,R_{9}$ in Figure~\ref{figAerial}.  Our first task is to
specify the prior distribution for $\mytheta$ and $\q$.  Following
\eqref{eqnPrior}, we take the prior for $p, z_0, L$ to be independent
of the prior on $\q$.  Since we do not have any prior information
about the values of the parameters $\pars$ (other than an acceptable
range), we choose a uniform distribution for each parameter
\begin{gather*}
  \prior(p)    \propto \mathbf{1}_{[0,0.6]}(p),\qquad
  \prior(z_{0})\propto \mathbf{1}_{[0,3]}(z_0),\qquad
  \prior(L)    \propto \mathbf{1}_{[-600,0]}(L),
\end{gather*}
where $\mathbf{1}_{E}$ is the indicator function for the interval set
$E$.  For the source strengths $\q$, we have the engineering estimates
from previous independent studies shown in Table~\ref{tabSources} which
can be applied along with a non-negativity constraint to construct a
suitable prior. In particular, we propose the following prior for source
strengths based on a gamma distribution:
\begin{equation}\label{gamma-prior-def}
  q_{j} \propto q_{j}^{\alpha_{j}-1}\exp(-\beta_{j}q_{j}).
\end{equation}
The quantities $\alpha_j$ are called shape parameters, while the
$\beta_j$ are rate parameters for the gamma distribution.  Values for
$\alpha_{j}$ and $\beta_{j}$ are derived from the engineering estimates
$\q=(q_{\textit{eng,j}})$ as the solution to the following system of
equations:
\begin{align*}
  \frac{\alpha_{j}-1}{\beta_{j}}             &=  q_{\textit{eng,j}}, \\
  \mathit{qgamma}(0.99,\alpha_{j},\beta_{j}) &= 3q_{\textit{eng,j}}.
\end{align*}
Here, $\mathit{qgamma}$ is the \emph{quantile function}, which is just
the inverse of the cumulative distribution function of the gamma
density.  Choosing $\alpha_{j}$ and $\beta_{j}$ in this way guarantees
that the mode of the prior distribution for $q_{j}$ is located at the
corresponding engineering estimate. Furthermore, $99\%$ of the prior
mass is concentrated in the interval $(0, 3q_{\textit{eng,j}})$ for each
$q_j$ \cite[Sec.~2.2]{gilchrist2000statistical}.

\begin{table}[htbp]
  \centering
  \caption{Engineering estimates of the source strengths from previous
    independent studies.}
  \label{tabSources}
  \begin{tabular}{|c|c|}
    \hline 
    Source & Estimated emission rate {[}ton/yr{]}\tabularnewline
    \hline\hline 
    $q_{\textit{eng},1}$ & 35\tabularnewline
    \hline 
    $q_{\textit{eng},2}$ & 80\tabularnewline
    \hline 
    $q_{\textit{eng},3}$ & 5\tabularnewline
    \hline 
    $q_{\textit{eng},4}$ & 5\tabularnewline
    \hline 
  \end{tabular}
\end{table}

Our next task is to specify the covariance matrix $\Sigma$ for the
measurement noise in \eqref{eqnNoisy}. Because the dustfall jars are
well-separated, we assume that the measurement noise for each jar is
independent from the others and hence take
\begin{gather}
  \label{eqnwhereLambda}
  \Sigma = \lambda I_{9\times 9}, 
\end{gather}
where $\lambda > 0$ is the variance of the measurement noise and
$I_{9\times 9}$ is the $9 \times 9$ identity matrix.  The variance
parameter $\lambda$ is related to the signal-to-noise ratio (SNR) of the
measurements and can be estimated by minimizing the functional
\begin{gather}
  \label{eqnFunctional2Optimize}
  J(\lambda) = \frac{1}{2} \int\left(
    \left\|\widehat{A}(p,L,z_{0})\q-\w) \right\|_{2} +
    \left\|\q-\q_{\textit{eng}}\right\|_{2}\right)  
  \,d\post^{\lambda},
\end{gather}
where $\q_{\textit{eng}}=[35,80,5,5]^{T}\,{ton}/{yr}$, and the notation
$\post^{\lambda}$ is used to explicitly show the dependence of the
posterior measure $\post$ on the parameter $\lambda$.  The motivation
for defining $J(\lambda)$ in this way comes from the problem of choosing
$\lambda$ to minimize the average value of the expression
\begin{gather}
  \label{eqnPreFunctional}
  (1-\delta) \left\|A(p,L,z_{0})\q-\w)\right\|_{2} + \delta
  \left\|\q-\q_{\textit{eng}}\right\|_{2} \qquad\text{for }\delta\in [0,1].
\end{gather}
Depending on the value of $\delta$, different weight can be given to the
credibility of the atmospheric dispersion model compared to that of
prior information about $\q$. In particular, if $\delta=1$ then
\eqref{eqnPreFunctional} reduces to
$\left\|\q-\q_{\textit{eng}}\right\|_{2}$ and $\lambda$ is chosen to
solely match the engineering estimates; on the other hand, taking
$\delta=0$ leaves only the term $\left\|A(p,L,z_{0})\q-\w)\right\|_{2}$
which chooses $\lambda$ solely based on minimizing the data misfit.  We
choose $\delta=\frac{1}{2}$, which is a compromise between the two
extremes.  After substituting into~\eqref{eqnPreFunctional} and
minimizing the expected value with respect to the posterior measure, 
the expression $J(\lambda)$ from~\eqref{eqnFunctional2Optimize} is
easily obtained.

\begin{figure}
  \centering
  \includegraphics[width=0.5\textwidth]{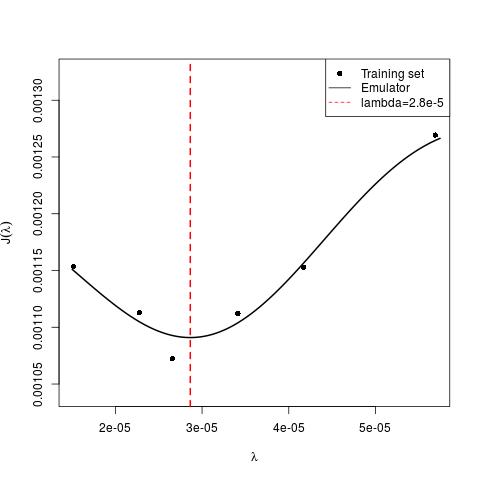}
  \caption{Emulator for $J(\lambda)$ in
    equation~\eqref{eqnFunctional2Optimize}.}
  \label{figLambdaEmul}
\end{figure}

Computing $J(\lambda)$ analytically is not practically feasible and so
we approximate it using a GP emulator, minimizing the emulator instead.
To obtain design points for the GP, we evaluate $J$ for six different
values of $\lambda$ using Markov Chain Monte Carlo integration, choosing
only six points to keep the computational cost within reasonable
limits. The design points were chosen empirically by identifying
locations where the function $J$ changes the most, and the resulting
simulation is shown in Figure~\ref{figLambdaEmul}. The minimum value of
$J(\lambda)$ is attained at $\lambda \approx 2.8\times 10^{-5}$, which
is the value of $\lambda$ we use in our source inversion framework
for~\eqref{eqnwhereLambda}.  Note that optimization with emulators is an
emerging area of research for which a detailed background can be found
in the paper \cite{osborne2009gaussian} and references therein.

We now apply Algorithm~\ref{algMH} to generate samples from the
posterior distribution, using the GP emulator $A(\pars)$ to compute the
likelihood ratio in Step 10 and parameter values listed in
Table~\ref{tab:MCMC-parameters}. When estimating integrals, we discard
the first half of the samples (up to 500,000) and retain only every
tenth iterate from there onward to reduce correlation between
consecutive samples. Therefore, a total of 50,000 samples is used to
approximate each integral with respect to the posterior.

\begin{figure}[htbp]
  \centering
    \includegraphics[width=0.49\textwidth]{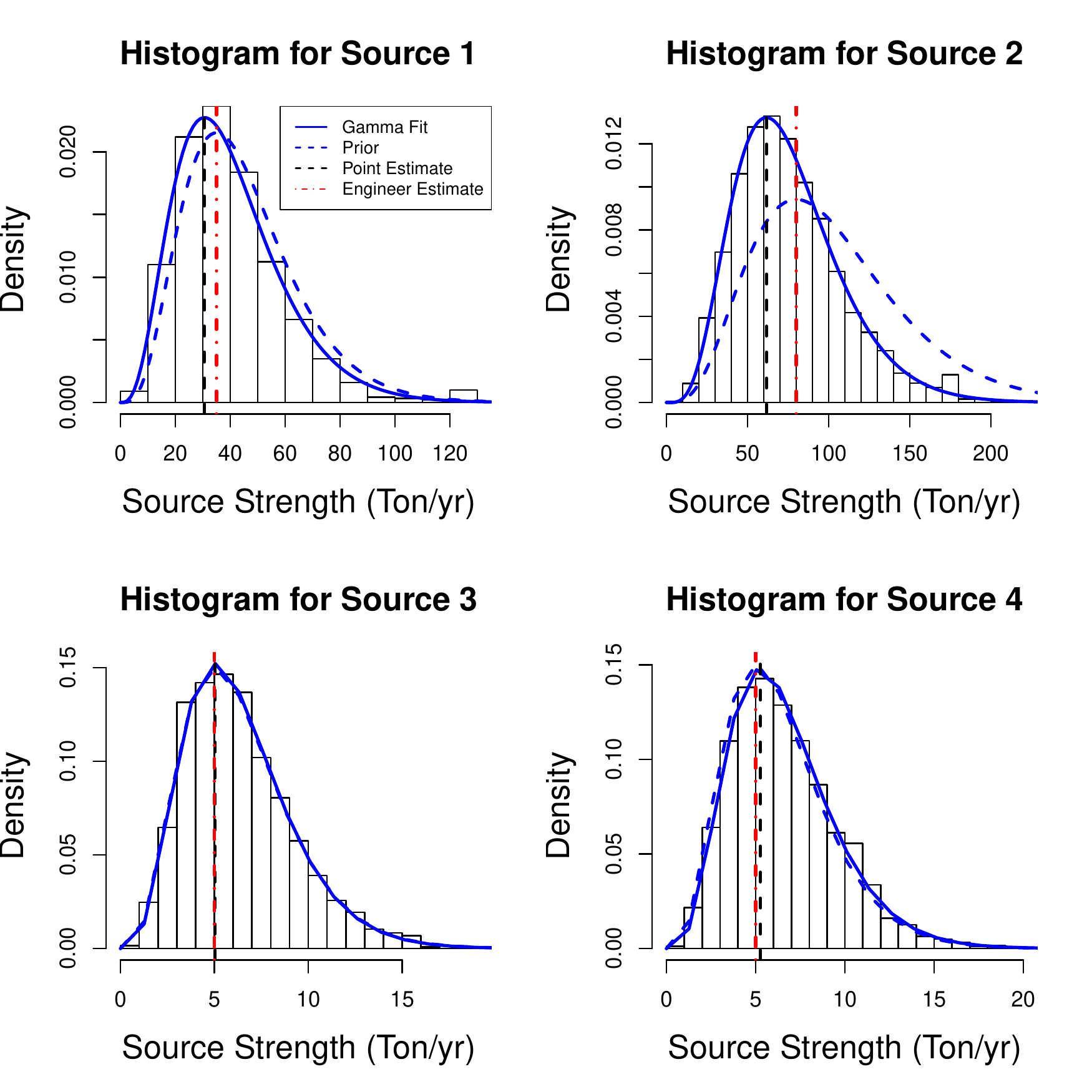}
  \includegraphics[width=0.49\textwidth]{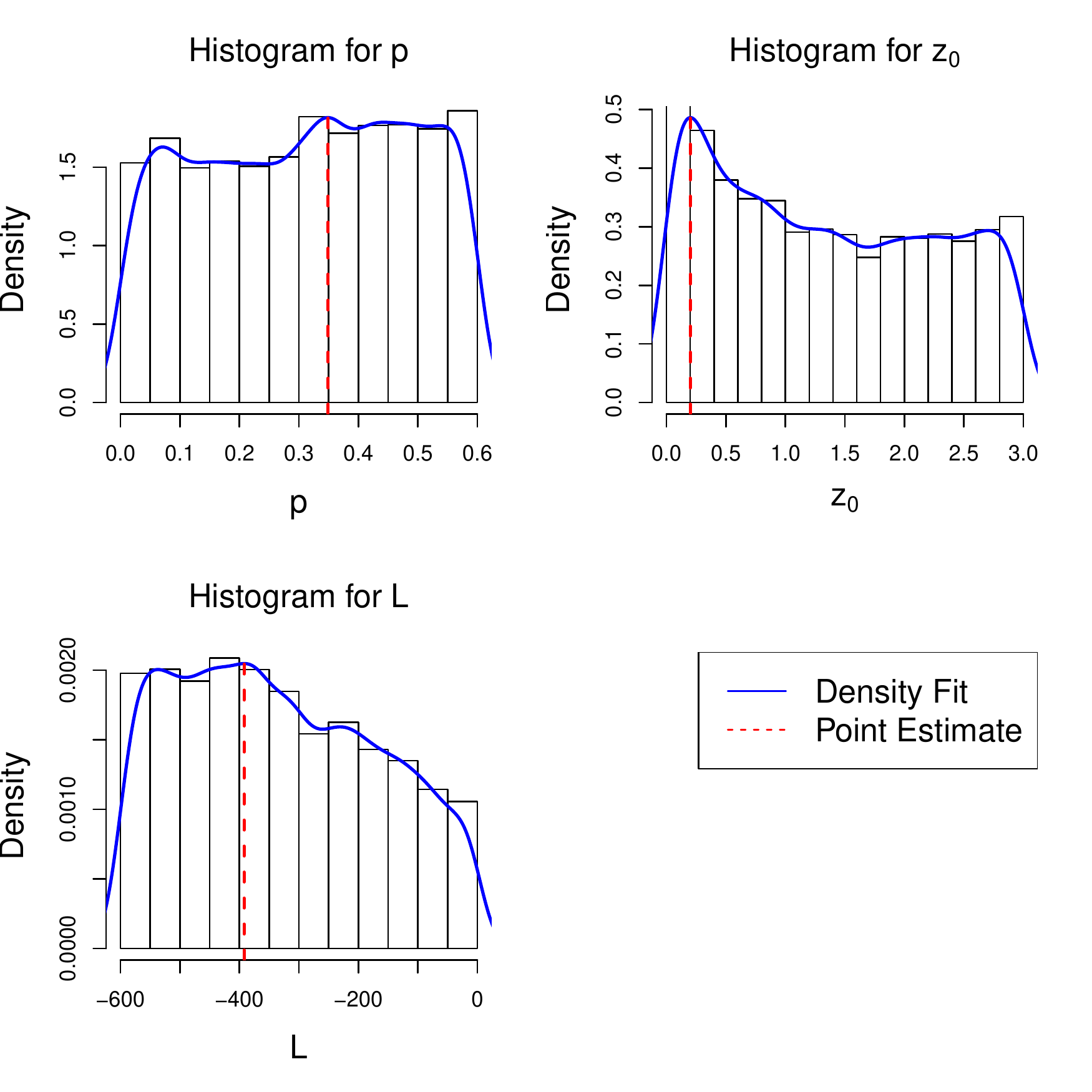}
  \caption{Marginal posterior distribution for the three parameters $\pars$
and the  four emission rates
    $q_{1},\dots, q_{4}$.}
  \label{figIndustrialCase}
\end{figure}

Marginal posterior histograms are displayed in
Figure~\ref{figIndustrialCase}, from which we observe that marginals for
the emission rates $q_j$ are skewed and have well-defined modes.  This
is in line with our choice of a gamma prior, and the mode for these
marginals is a reasonable point estimate for each $q_j$. To approximate
the marginal mode we fit a gamma distribution to the histograms and take
the point estimate to be the mode of the fitted distribution.

\begin{table}[htbp]
  \centering
  \caption{Value of parameters used in Algorithm~\ref{algMH} for
    generating samples from the posterior distribution.} 
  \label{tab:MCMC-parameters}
  \begin{tabular}{|c|c|c|c|c|}
    \hline
    $N$    & $\beta$ & $\gamma_1$& $\gamma_2$ & $\gamma_3$ \\ \hline
    $10^6$ & $0.05$  & $0.01$    & $(2.38)^2$ & $(0.1)^2$ \\ \hline
  \end{tabular}
\end{table}

The situation is slightly different for the parameters $\pars$.  The
histogram for $p$ has no distinctive maximum point indicating that the
data is not informative for this parameter. For $z_{0}$ and $L$, some
values are more distinctive than others but the distinction is still not
sharp. To choose a point estimate for these cases, we first fit a
density over the histograms and then take the resulting mode as the point
estimate, which is indicated in each case by a red dashed line in
Figure~\ref{figIndustrialCase}.

\begin{figure}[htbp]
  \centering
  \includegraphics[width=0.5\textwidth]{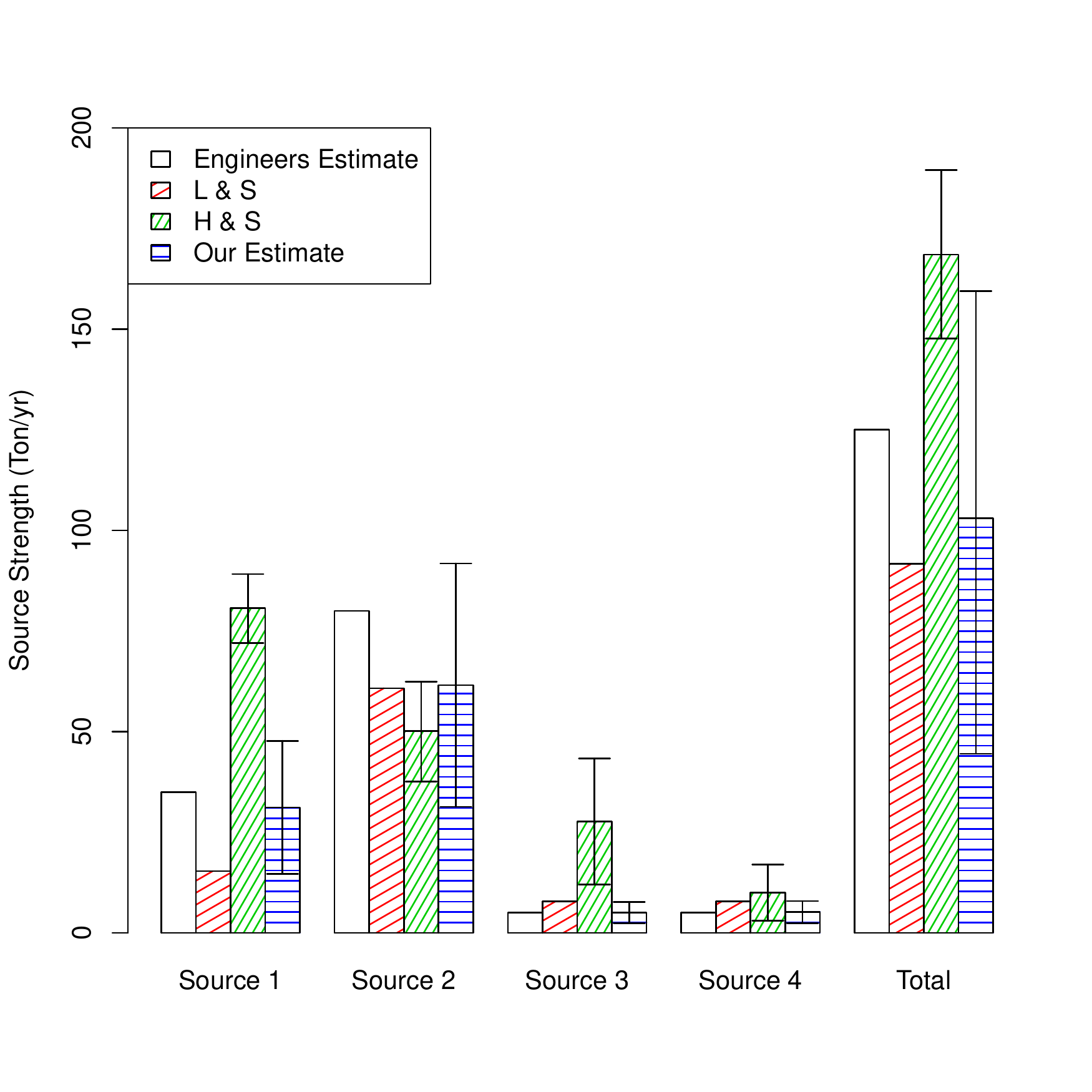}
  \caption{Comparison between engineering estimates of the emission
    rates and our point estimates, along with the associated
    uncertainty. The ``Total'' column indicates the sum of emission
    rates from all four sources. Previous results from Lushi and
    Stockie~\cite{lushi2010inverse} (L\&S) and Hosseini and
    Stockie~\cite{hosseini2016airborne} (H\&S) are included for
    comparison purposes.}
  \label{figQUnicertainty}
\end{figure}

To quantify the uncertainty in each parameter, we employ a $68\%$
Bayesian confidence interval which is defined as follows: given a point
estimate $x^{\ast}$ of a random variable $X$ distributed with
probability density $\rho$, a $68\%$ Bayesian confidence interval is the
radius of the ball centered at $x^{\ast}$ containing $68\%$ of the
probability mass for $\rho$.  The corresponding confidence interval is
listed for each parameter in Table~\ref{tabFinalEstimates}.  The
estimates of the $q_j$ are clearly more informative, which is due to the
linearity of the forward map and stronger prior knowledge. The marginal
modes and associated uncertainties for the $q_j$ are depicted
graphically in Figure~\ref{figQUnicertainty}.
These results suggest that the engineering estimates for sources $q_1$
and $q_2$ are most likely overestimates, whereas sources $q_3$ and $q_4$
are slightly underestimated.  Furthermore, our estimates for sources
$q_{2}$, $q_{3}$ and $q_{4}$ qualitatively agree with the results of
previous studies in \cite{lushi2010inverse} (L\&S) and
\cite{hosseini2016airborne} (H\&S).

The results in (L\&S) were obtained using a Gaussian plume
model and a least-squares based inversion method. Therefore,
discrepancies between our estimates and those of (L\&S) are expected
given that we used a very different forward model.  On the other hand, (H\&S)
used the same PDE based forward map as this study (although within a
different Bayesian formulation) and so it is more informative to compare
our results to those of (H\&S).  We note two key differences between our
results.  First of all, there is a notable difference between the
estimates for $q_1$.  In \cite{hosseini2016airborne}, the parameters
$p,z_{0},L$ are fixed at $0.2, 0.05$ and $-10$, respectively.  But these
values depend strongly on environmental conditions, so a mis-assessment
in the environmental conditions might lead to sub-optimal
predictions. It was conjectured in \cite{hosseini2016airborne} that this
may be the underlying reason for the apparent overestimation of $q_1$ in
comparison to previous studies. Furthermore, (H\&S) obtained their
estimates by marginalizing the posterior measure on $q_j$ only, whereas
we marginalize the posterior measure on the $q_j$ \emph{and} the model
parameters $\mytheta$. \hl{We believe that the discrepancy between our
  results and \cite{hosseini2016airborne} highlights the importance
  of model calibration in atmospheric source inversion applications.}


The second major difference is that our estimates have larger error bars
in comparison to (H\&S), although we note that our error bars are large
enough to cover most previous estimates.  This is largely due to
differences in the parameter $\lambda$ which have a direct effect on
posterior variance.  We chose $\lambda$ by minimizing $J(\lambda)$ in
\eqref{eqnFunctional2Optimize} and achieved an SNR of $2.70$. This means
that the strength of the signal is quite low compared to the noise of
the measurements. In \cite{hosseini2016airborne}, $\lambda$ is chosen to
achieve an SNR of 10 which is why the error bars for
those results are much smaller. Furthermore, our method accounts for
uncertainties in the value of the model parameters $\mytheta$ while
these are fixed in the approach of (H\&S).

We further observe that our estimate of the total emission is consistent
with previous studies, with the exception of (H\&S) which appears to
have a higher estimate. This is mostly due to the over-estimation of
$q_1$ which, as was mentioned before, is likely due to a lack of
calibration in \cite{hosseini2016airborne}.  By allowing the model
parameters to be calibrated automatically we have obtained estimates
that are consistent with two previous studies that used very different
methodologies, namely the engineering estimates and (L\&S). This
highlights the importance of model calibration in atmospheric source
inversion studies.  Furthermore, our uncertainty estimates are more
realistic since we account for model parameter uncertainties and infer a
better value of the signal-to-noise ratio.  Although our uncertainty
estimate for the total emission is relatively large in comparison to the
estimate itself, this is due to the fact atmospheric source inversion
problems are severely ill-posed and sensitive to model parameters.

\begin{table}[htbp]
  \centering
  \caption{Parameters and their estimates and uncertainties.}
  \label{tabFinalEstimates}
  \begin{tabular}{|c|c|c|c|}
    \hline 
    Parameter & Point Estimate & $68\%$ Confidence Interval\tabularnewline
    \hline 
    \hline 
    $p$ &  0.3478 & $[0.1498,0.5458]$\tabularnewline
    \hline 
    $z_{0}$ & 0.0811 & $[0,1.5781]$\tabularnewline
    \hline 
    $L$  & $-$379.45 & $[-195.86,-563.04]$\tabularnewline
    \hline 
  \end{tabular}
\end{table}

\hll{Finally, we note that GP emulation is a crucial
  aspect of our inversion method that renders it cost-effective. The
  main computational bottleneck in our framework is the evaluations of
  the finite volume forward solver at the 64 design points used to
  construct the emulator. Given that the design points are independent
  not only of each other but also the measurements, these computations
  can be performed offline and in parallel.  Every run of the finite
  volume code requires roughly two hours on a personal computer, while
  the emulator can be evaluated in a fraction of a second.  Therefore,
  the cost of evaluating the emulator is negligible in comparison to the
  finite volume solver which enables us to perform Monte Carlo
  integration in a reasonable time. Clearly, such a calculation would be
  prohibitive had we employed the finite volume solver directly in the
  MCMC step.

  In comparison to the previous approaches of
  \cite{hosseini2016airborne} and \cite{lushi2010inverse} our total
  computational cost (including the training of the emulator and design
  of experiment) is higher. For example, \cite{lushi2010inverse} uses a
  Gaussian plume model which is much cheaper than the finite volume
  solver but comes at the cost of major simplifications to the physics.
  On the other hand, \cite{hosseini2016airborne} uses a finite volume
  forward solver of the original PDEs, which is only evaluated four
  times. However, neither of these two approaches provides an automated
  calibration procedure, and so we claim that the extra computations in
  our methodology are a small price to pay for a richer solution
  structure with more realistic uncertainty estimates and rigorous model
  calibration.  }

\subsection{Effect of the prior choice}\label{sec:effect-prior-choice}

\hl{Here we study how the choice of the prior on emission rates affects
  the marginal posteriors.  We use the same uniform priors on the model
  parameters $\pars$ as before because those uniform priors are readily
  uninformative. We take the gamma priors of
  Section~\ref{sec:source-inversion} and vary the parameters $\alpha_j$
  and $\beta_j$ in \eqref{gamma-prior-def}.  More precisely, we choose
  $\alpha_j, \beta_j$ to solve
  \begin{align*}
    \frac{\alpha_{j}-1}{\beta_{j}}             &=  q_{\textit{eng,j}}, \\
    \mathit{qgamma}(0.99,\alpha_{j},\beta_{j}) &= \tau q_{\textit{eng,j}}.
  \end{align*}
  with $\tau = 2, 3, 4$ and $qgamma$ being the quantile function for the
  Gamma distribution as before. Recall that these equations imply that
  all prior modes are centered at the engineering estimates but $\tau$
  controls the spread of the priors.  In fact, 99\% of the prior mass
  lies in the interval $(0, \tau q_{\textit{eng},j})$.  Then setting
  $\tau =2$ forces the prior to concentrate more tightly around the mode
  while $\tau=4$ corresponds to a flat
  prior. Figure~\ref{figPriorEffect} shows the three choices of the
  priors for each source as well as the corresponding marginal
  posteriors. For sources 1 and 2 a noticeable contraction of the
  posterior occurs, which is more noticeable for larger $\tau$. The
  contraction of the posterior implies that the data is in fact
  informative regarding the first two sources.  On the other hand, the
  marginals on sources 3 and 4 are nearly identical to their priors,
  which is a sign that the data is not informative about these
  sources. We note that this observation is not surprising given the
  fact that the emissions from sources 3 and 4 are negligible in
  comparison to 1 and 2. 

  In Figure~\ref{figPriorEffectEmissionRates} we compare the pointwise
  estimators of the emission rates as well as uncertainty estimates for
  different choices of $\tau$.  As expected, the uncertainties increase
  with $\tau$ since prior uncertainty directly affects posterior
  uncertainty. We also observe that smaller values of $\tau$ pull the
  pointwise estimates towards the engineering estimates, which is also
  in line with the intuition that choosing a prior that is too tight can
  result in the posterior being overwhelmed by the prior which tends to
  ignore the measured data.  Furthermore, by comparing to
  Figure~\ref{figPriorEffect} we observe stronger posterior contraction
  relative to prior for the first two sources in the case where $\tau$
  is larger. This indicates that the measurements are indeed informative
  in the direction of the larger sources. Relative contraction of the
  posterior appears to be weaker when $\tau =2$ indicating that the
  tighter prior is too strong.

  \begin{figure}[htbp]
    \centering
    \includegraphics[width=0.49\textwidth]{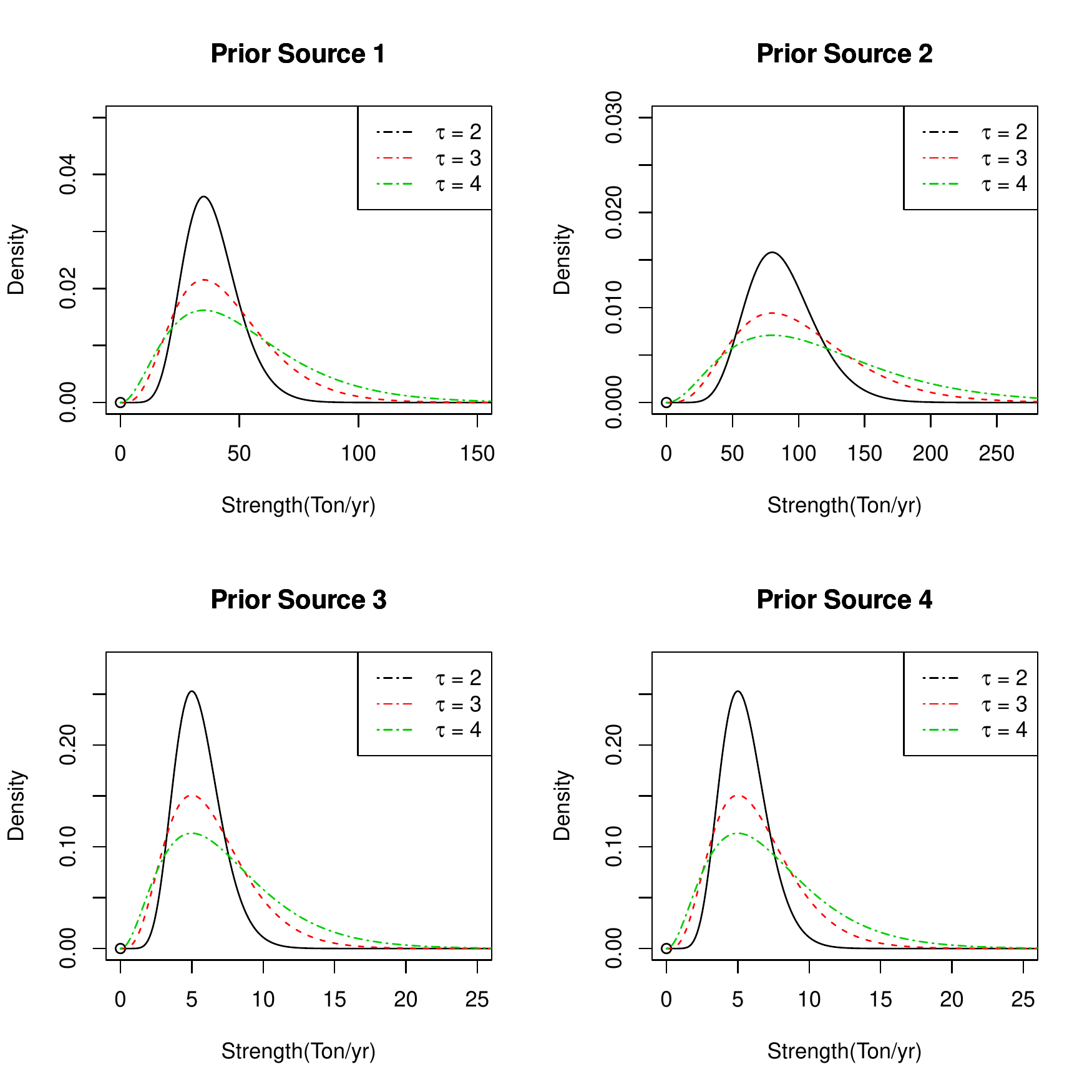}
    \includegraphics[width=0.49\textwidth]{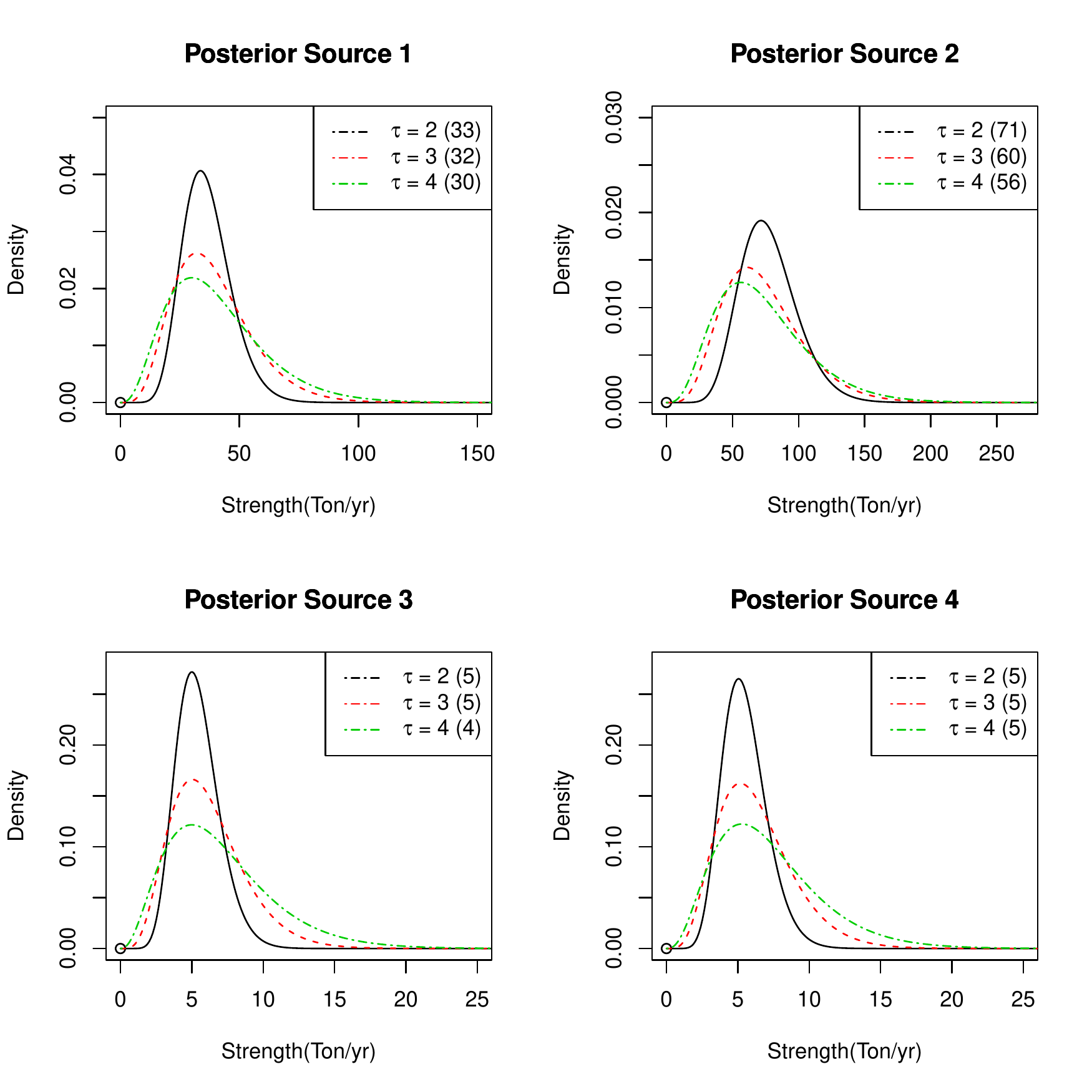}
    \caption{Comparison between different choices of prior
      distribution on the emission rates and the corresponding
      posterior distributions for $\tau =2, 3, 4$. Smaller values of
      $\tau$ give tighter priors that are more concentrated around the
      engineering estimates while larger values result in more flat
      priors that are less informative.}
    \label{figPriorEffect}
  \end{figure}
  
  \begin{figure}[htbp]
    \centering
    \includegraphics[width=0.5\textwidth]{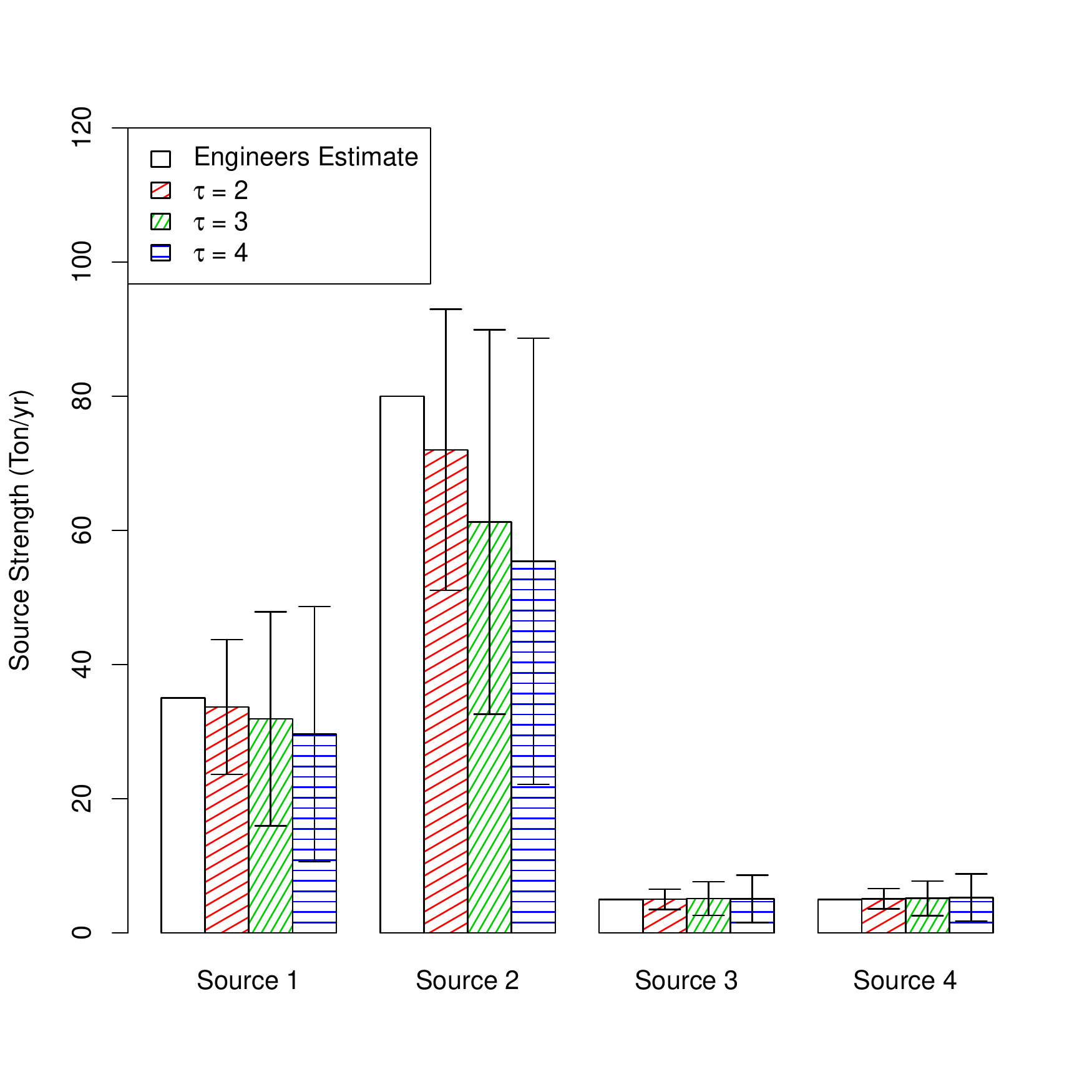}
    \caption{Comparison between point estimates of emission rates 
      using different prior distributions, along with the associated
      uncertainties.}
    \label{figPriorEffectEmissionRates}
  \end{figure}
}

\subsection{Effect of emulator quality}\label{sec:effect-emul-qual}

\begin{figure}[htbp]
  \centering
  \includegraphics[width=0.49\textwidth]{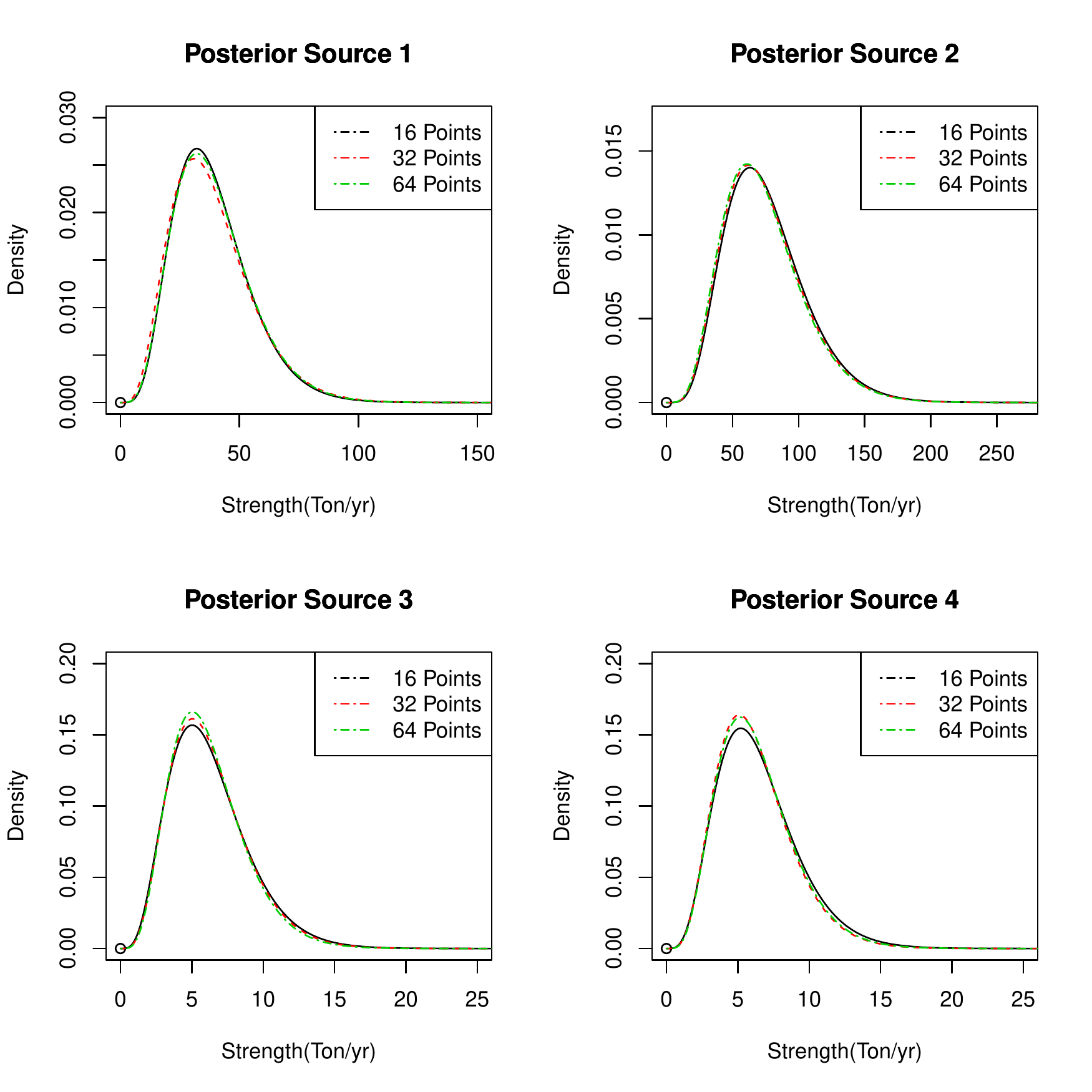}
  \includegraphics[width=0.49\textwidth]{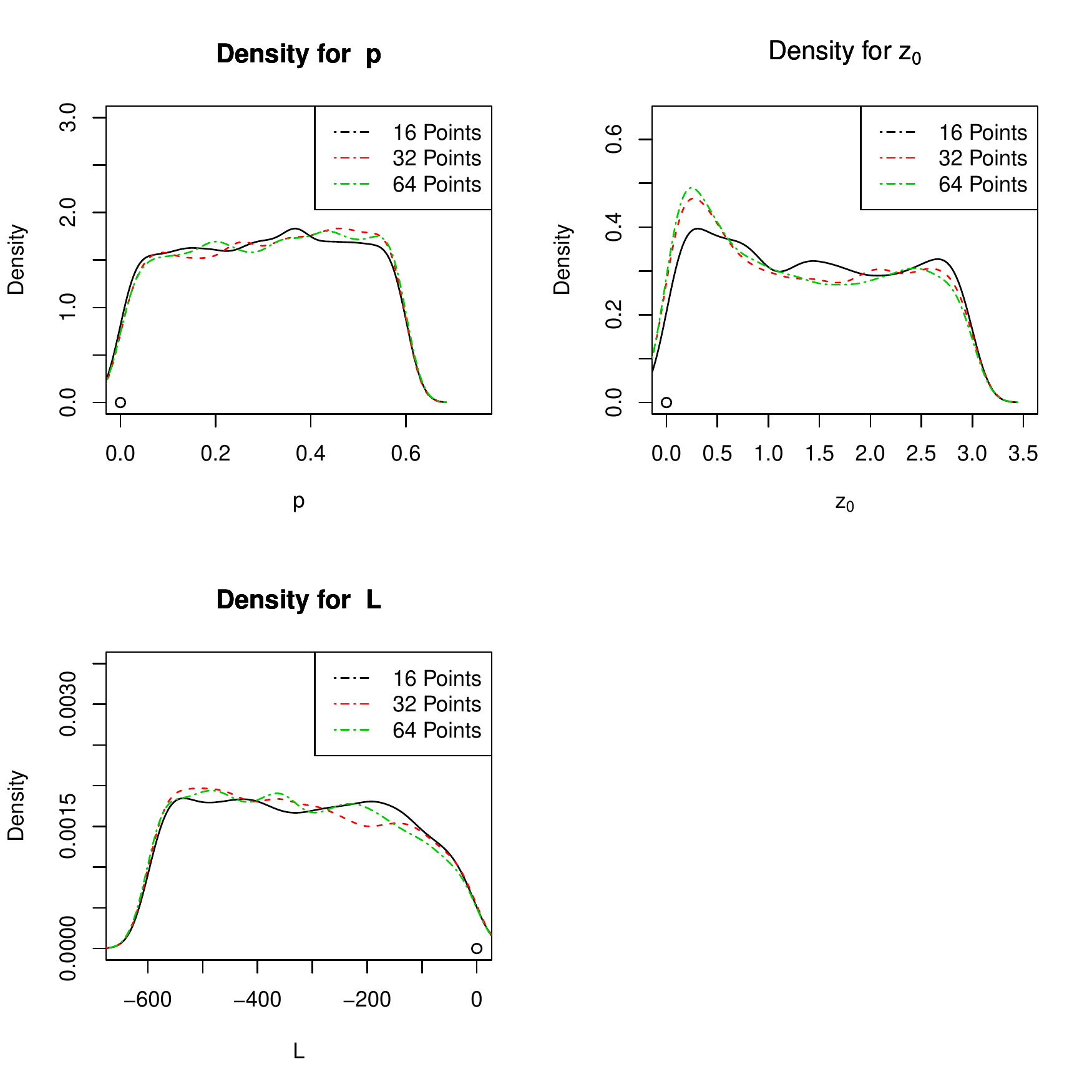}
  \caption{Marginal posterior distributions of model parameters $p, L,
    z_0$ and the four emission rates $q_1, \dots, q_4$, for different
    emulators using experimental designs with 16, 32 and 64 points.}
  \label{figEmulatorConvergence}
\end{figure}

\hl{In this section we demonstrate how the quality of the GP emulator
  affects the posterior distributions in
  Figure~\ref{figIndustrialCase}. We construct a hierarchy of GP
  emulators using 16, 32 and 64 experimental designs and then solve the
  inverse problem with the corresponding emulators. Our results are
  summarized in Figure~\ref{figEmulatorConvergence} where we plot the
  density fit to the posterior marginals of the model parameters and
  emission rates.  We observe that the marginals on the emission rates
  have essentially converged and there is very little difference between
  the different emulators.  More variations are present between the
  marginals on the model parameters $p, L$ and $z_0$.  In the case of
  $p$ and $L$ the small perturbations are most likely due to numerical
  errors in computing the histograms and fitting the densities. In the
  case of $z_0$ we see larger variations, specifically in the
  maximal point emerging close to 0.25 which is present when 32 and 64
  point experimental designs are utilized. The closeness of the 32 and
  64 point density fits indicate convergence of the posterior marginals
  for $z_0$ and provide further confidence that the maximal point is
  indeed located near $0.25$.}

\section{Conclusions}

We have developed a new method for simultaneously calibrating
atmospheric dispersion models and solving source inversion
problems. GP emulators are used to speed up the evaluation
of the atmospheric dispersion model, which are then employed within a
Bayesian framework to calibrate model parameters and estimate source
strengths at the same time. This allows for automatic tuning of unknown
model parameters rather than more heuristic approaches that are
typically used in practice.

We demonstrated the effectiveness of our proposed method in an
industrial case study concerning emissions of zinc from a lead-zinc
smelter in Trail, BC, Canada.  Our results agree with those of previous
studies in the literature that are based on very different
approaches. Hence, we conclude that the emulation process is an accurate
and efficient one that has the added advantage of greatly speeding up
the calculations by exploiting emulation.

Our proposed solution methodology is not restricted to atmospheric
dispersion models governed by partial differential equations. Indeed,
our method can be adapted to any application where one in dealing with a
forward map that is linear in the primarily unknowns (such as the source
strengths) but may depend nonlinearly on the model parameters.

\begin{acknowledgements}
  This work was partially supported by the Natural Sciences and
  Engineering Research Council of Canada through a Postdoctoral
  Fellowship (BH) and a Discovery Grant (JMS).  We are grateful to the
  Environmental Management Group at Teck Resources Ltd.\ (Trail, BC) for
  providing data and for many useful discussions.
\end{acknowledgements}

\bibliography{paagPaperV8}
\bibliographystyle{abbrv}

\end{document}